\documentclass[twoside,11pt]{article}
\usepackage[top=1in, bottom=1in, left=1in, right=1in]{geometry}
 \usepackage[backref,colorlinks,linkcolor=red,anchorcolor=green,citecolor=blue]{hyperref}
\usepackage{amsfonts,amssymb}
\usepackage{amsmath}
\usepackage{graphicx}
\usepackage{cite}
\usepackage{enumerate}
\usepackage{caption}
\usepackage{subcaption}
\usepackage{booktabs} 
\usepackage{multirow}
\usepackage[ruled]{algorithm2e}
\usepackage{fancyhdr}

\renewcommand{\theequation}{\arabic{section}.\arabic{equation}}

\numberwithin{equation}{section}

\newcommand{\In}{i}
\newcommand{\Jn}{j}

\newcommand{\I}{\mathcal{I}}
\newcommand{\J}{\mathcal{J}}
\newcommand{\Ir}{\widetilde{\I}}
\newcommand{\Is}{\hat{\I}}

\renewcommand{\S}{\mathcal{S}}
\newcommand{\Bi}{\mathcal{B}}
\newcommand{\R}{\mathcal{R}}
\newcommand{\Rb}{\mathbb{R}}

\newcommand{\A}{{A}}
\newcommand{\Hm}{{A}}
\newcommand{\Hs}{{\widetilde{A}}}
\newcommand{\He}{{\bar{A}}}
\renewcommand{\AE}{\widetilde{\A}}
\newcommand{\E}{{M}}
\newcommand{\F}{\mathcal{F}}
\newcommand{\T}{T}
\newcommand{\Tree}{\mathcal{T}}

\newcommand{\K}{{K}}
\newcommand{\Kt}{K^{*}}
\renewcommand{\L}{{L}}

\newcommand{\U}{{U}}
\newcommand{\Ub}{\widetilde{U}}
\newcommand{\B}{{B}}

\newcommand{\adjkernel}{S}

\newcommand{\Z}{{Z}}
\newcommand{\C}{{C}}

\newcommand{\Ka}{\dot{K}}
\newcommand{\Sa}{\dot{S}}
\newcommand{\Ga}{\dot{G}}

\def\@begintheorem#1#2{\par\bgroup{\noindent\scshape #1\ #2. }\it\ignorespaces}
\def\@opargbegintheorem#1#2#3{\par\bgroup%
   {\scshape #1\ #2\ ({\upshape #3}). }\it\ignorespaces}
\def\@endtheorem{\egroup}

\newtheorem{theorem}{Theorem}
\newtheorem{proposition}[theorem]{Proposition}

\newlength{\bibitemsep}
\newlength{\bibparskip}
\let\oldthebibliography\thebibliography
\renewcommand\thebibliography[1]{%
  \oldthebibliography{#1}%
  \setlength{\parskip}{\bibitemsep}%
  \setlength{\itemsep}{\bibparskip}%
  \setlength\itemindent{-15pt}
}

\allowdisplaybreaks
\begin{document}
\title{Approximate inversion of discrete Fourier integral operators}


\author{Jordi Feliu-Fab\`{a}\thanks{Institute for
    Computational and Mathematical Engineering, Stanford University,
    Stanford, CA 94305, (jfeliu@stanford.edu).}
  \and
  Lexing Ying\thanks{Department of Mathematics and Institute for Computational and Mathematical
    Engineering, Stanford University, Stanford, CA 94305, (lexing@stanford.edu)}
}

\pagestyle{fancy}
\fancyhf[HLO]{Approximate inversion of discrete Fourier integral operators}\fancyhf[HRE]{Jordi Feliu-Fab\`{a}, Lexing Ying}
\fancyhf[HRO,HLE]{\thepage}

\date{}
\maketitle

\begin{abstract}
  This paper introduces a factorization for the inverse of discrete Fourier
  integral operators that can be applied in quasi-linear time. The factorization starts by
  approximating the operator with the butterfly factorization. Next, a hierarchical matrix
  representation is constructed for the hermitian matrix arising from composing the Fourier integral
  operator with its adjoint. This representation is inverted efficiently with a new algorithm based on the hierarchical interpolative factorization. By combining these two factorizations, an approximate inverse
  factorization for the Fourier integral operator is obtained as a product of $O(\log N)$ sparse
  matrices of size $N\times N$. The resulting approximate inverse factorization can be used as a
  direct solver or as a preconditioner. Numerical examples on 1D and 2D Fourier integral operators,
  including a generalized Radon transform, demonstrate the performance of this new approach.
\end{abstract}

{\bf Keywords:} Fourier integral operator; Radon transform; hierarchical matrices; hierarchical interpolative factorization; butterfly algorithm

{\bf AMS subject classifications:} 35S30; 44A12; 65F05

\section{Introduction}\label{introduction}

In this paper, we consider Fourier integral operators (FIOs) of the form
\begin{align}
  \label{eq1}
  (\mathcal{L} f)(x) = \int_{\Rb^d} a(x,\xi)e^{2\pi i \Phi(x,\xi)}\hat{f}(\xi)\text{d}\xi.
\end{align}
Here, $\xi$ is the frequency variable, $x$ the spatial variable, $a(x,\xi)$ is a smooth amplitude function in both $x$ and $\xi$, $\Phi(x,\xi)$ is a smooth phase function in $(x,\xi)$ for $\xi\neq 0$ with the homogeneity relation  $\Phi(x,\lambda\xi) = \lambda\Phi(x,\xi)$ for $\lambda$ positive, and $\hat{f}$ is the Fourier transform of $f$ defined by
\begin{align}
  \hat{f}(\xi) = \int_{\Rb^d} e^{-2\pi i x\cdot\xi}f(x)\text{d}x.
\end{align}


Typical uses of FIOs include the study of hyperbolic equations or integration over certain curved
manifolds used in inverse problems. For instance, FIOs are commonly used to represent the solution
to wave equations in $\Rb^d$ \cite{Demanet}.

To introduce the discrete analog of an FIO in the $d$-dimensional case, we start with two Cartesian grids
\begin{align}
  & X =\left\{x = \left(\frac{n_1}{n} , \ldots, \frac{n_d}{n} \right),0 \leq n_1,\ldots, n_d < n\text{ and }n_1, \ldots n_d \in \mathbb{Z}\right\}, 
  \label{eq:X_domain}\\
  & \Omega = \left\{ \xi=(n_1,\ldots, n_d),-n/2 \leq n_1,\ldots,n_d < n/2 \text{ with } n_1,\ldots,n_d \in \mathbb{Z}\right\},
    \label{eq:freq_domain}
\end{align}
where the number of discretization points per dimension $n$ is assumed to be even for simplicity.  The $d$-dimensional discrete Fourier integral operator of a function $f(x)$ is then defined by
\begin{align}
  \label{eq3}
  (\mathcal{L} f)(x) \equiv (\mathcal{K}\hat{f})(x) =  \sum_{\xi\in \Omega} a(x,\xi)e^{2\pi i \Phi(x,\xi)}\hat{f}(\xi)
\end{align} 
for every $x\in\mathcal{X}$, where $\hat{f}$ is the discrete Fourier transform of $f$ defined by
\begin{align}
  \hat{f}(\xi) = \frac{1}{n^d}\sum_{x\in X} e^{-2\pi i x\cdot\xi}f(x).
\end{align}


The discrete FIO from \eqref{eq3} can be written in matrix form as
\begin{align}
  \label{eq6}
  u = \K \hat{f}
\end{align}
with $\hat{f}=Ff$ and $u$ being the discrete Fourier transform and the discrete Fourier integral
transform of $f$, respectively, evaluated at each discretization point. Notice that $\K$ is a dense
matrix of size $n^d\times n^d$, with entries $\K_{ij} = a(x_i,\xi_j)e^{2\pi i \Phi(x_i,\xi_j)}$.

While in some cases, one is only interested on applying a discrete FIO quickly, in others one needs to solve a system of equations with a discrete FIO matrix, which in general is computationally challenging. In particular, FIOs are widely used in inverse problems, where given some measurements $u$ one wants to compute
some unknown quantity of interest $f$ by solving \eqref{eq6}. Some important applications, that can be formulated
as such, are radar imaging \cite{Yazici,Cheney}, thermoacoustic and photoacoustic tomography
\cite{Stefanov,Kuchment}, seismic imaging \cite{Symes}, single photon emission computed tomography
\cite{Kunyansky} and electrical impedance tomography \cite{santosa}, among many others.

Many matrix factorizations for the rapid application of discrete FIOs can be found in the literature
that reduce the computational complexity from quadratic $O(N^2)$ to quasi-linear $O(N\log N)$ in the
number of discretization points $N = n^d$. However, those factorizations are not readily
invertible. In this paper, we are concerned with providing a factorization for the inverse of discrete
FIOs that can be applied in quasi-linear time in order to solve \eqref{eq6} efficiently. One can
either solve $u = \K \hat{f}$ or  $u = \K Ff$ by  inverting $K$ or $KF$. In our numerical examples, we observe that
inverting $K$ is more efficient than inverting $KF$, since it leads to lower ranks and faster decay
of the entries during the factorization. Since $F$ is the discrete Fourier transform matrix, it can be
inverted and applied efficiently with fast algorithms such as the fast Fourier transform (FFT), therefore in the rest of the paper we will be concerned with
inverting $K$, which we refer to as FIO matrix.

The proposed factorization is tested on several FIOs, including a 2D generalized Radon Transform
(RT) \cite{Hu}. The Radon Transform can be seen as a set of line integrals, with many applications
in seismic data processing and various types of tomography, where the integrals can be taken along
lines, parabolas (parabolic RT), hyperbolas (hyperbolic RT) or spheres (spherical RT).

\subsection{Background.}\label{sec:background}

There has been significant amount of research on developing fast algorithms and factorizations for
the direct evaluation of discrete FIOs. By exploiting the complimentary low-rank condition of the
FIO matrix $\K$, these algorithms compress certain approximate low-rank matrices that appear
throughout the factorization process. Many of the state-of-the-art algorithms
\cite{BF1d,BF2d,BFId,IDBF,MIDBF,BFrand} can be viewed as factorization schemes based on the
butterfly algorithm \cite{BFalgo, Candes}.  Such factorizations, while efficient to represent $\K$,
cannot be easily inverted due to the telescoping nature of the factorization which leads to
rectangular matrices throughout the factorization. Other approaches include fast directional
algorithms \cite{Candes_dir}. While many algorithms have been proposed for the efficient
factorization of $\K$, we are not aware of any prior literature on the efficient factorization for
the inverse operator $\K^{-1}$ for a general class of FIOs.

A first approach to invert $\K$ consists on using classical direct methods such as Gaussian
elimination or other standard matrix factorizations \cite{Golub}, which unfortunately achieve
$O(N^{3})$ complexity.  For some dense matrices, this can be accelerated by exploiting the
low-rankness of certain submatrices. Many representations of dense matrices exist such as wavelet
decompositions \cite{Beylkin} or hierarchical matrices based on different admissibility conditions
and the use of nested bases, including $\mathcal{H}$-matrices \cite{Hackbusch,Hackbusch_Born,IFMM},
hierarchically semiseparable (HSS) matrices \cite{Xing,Aminfar}, and hierarchically off-diagonal
low-rank (HODLR) matrices \cite{Aminfar2}, among others \cite{Bebendorf:2013,Martinsson:2005,Minden,
  Corona}. These algorithms can be much faster than classical direct methods, with some even
attaining quasi-linear complexity to compute the inverse $\K^{-1}$. However, some of these
approaches rely on the effective compression of the operator with wavelets, while others rely on a
hierarchical partitioning of the domain and the admissibility condition at each level, which FIO
matrices typically don't satisfy with such partitioning.

Alternatively, one can solve the system of equations \eqref{eq6} using iterative methods
\cite{Saad}. That is with algorithms such as the generalized minimal residual method
(GMRES)\cite{Saad_GMRES} or the biconjugate gradient method (BiCG) \cite{Vorst}, which can achieve
$O(N^2)$ complexity if only a small number of iterations are required. The complexity can be further
reduced to $O(N\log N)$ by using the butterfly factorization to perform fast matrix-vector
multiplications with $\K$. However, in general, the number of iterations required by an iterative
solver is sensitive to the spectral properties of the matrix of the system to be solved. In the case
of FIO matrices, the number of iterations can be proportional to the problem size $N$, requiring
efficient preconditioners. For the FIO matrices used in the numerical tests of this
paper, the eigenvalues are located around a circle in the complex plane, and therefore, due to the
lack of eigenvalue clusters, convergence is slow. A particular preconditioner that can be used is
the adjoint FIO matrix $\Kt$, which reduces the number of iterations, but still converges slowly in some cases.

As a result of these observations, there is a need for finding algorithms that can compute an
approximate inverse factorization of FIO matrices or provide efficient preconditioners. The main
purpose of this paper is to address this issue.

\subsection{Contributions.}\label{sec:contrib}
The main contribution of this paper is the introduction of an algorithm to obtain a factorization
for the inverse FIO matrix $\K^{-1}$ as a product of $O(\log N )$ sparse matrices, each with $O(N)$
nonzero entries, that can be applied in quasi-linear time to a random vector.  The main
building blocks to build such a factorization are:
\begin{itemize}
\item Building the butterfly factorization approximation $\Ka$ of $K$ using the algorithm from
  \cite{BF1d,BF2d}.
\item Constructing the $\mathcal{H}$-matrix approximation $\Sa$ of $S\equiv K^*K$ using the peeling
  algorithm \cite{Peeling} and quasi-linear time matrix-vector multiplication with $\Ka$.
\item The $\mathcal{H}$-matrix approximation $\Sa$ is inverted using a new algorithm based on the
  hierarchical interpolative factorization \cite{HIFIE} proposed in this paper to obtain the
  approximate factorization $\Ga$ of $G \equiv (K^* K)^{-1}$.
\end{itemize}
The inverse factorization is then obtained as 
\[\K^{-1} \approx \Ga \Ka^{*}. \]
Throughout this paper, we shall use the dot accent to a matrix to indicate its approximation.

\subsection{Outline.}\label{sec:outline}
The rest of this paper is organized as follows. Sections \ref{sec:BF} and \ref{sec:Peeling} review
some of the tools used for constructing the inverse factorization, i.e. the butterfly factorization,
hierarchical matrices and the peeling algorithm. Section \ref{sec:HIF} presents the algorithm used
to invert hierarchical matrices. Section \ref{sec:Alg} proposes an algorithm for the inverse
factorization of FIOs by combining the previous factorizations from sections
\ref{sec:BF}, \ref{sec:Peeling} and \ref{sec:HIF}. In Section \ref{sec:NR}, numerical results are provided to demonstrate
the efficiency of the algorithm for 1D and 2D FIOs. Finally, Section \ref{sec:conclusions} concludes
with some discussion and future work.

\section{Butterfly factorization}\label{sec:BF}
In order to built the inverse approximation of an FIO matrix $\K$ we first need an algorithm to
perform matrix-vector multiplications with $\K$ and $\Kt$ efficiently. In this paper we use the
multidimensional butterfly factorization \cite{BF1d,BF2d} to obtain a factorization of the FIO
matrix $\K$ and its conjugate $\Kt$ that can be applied in quasi-linear time to a random
vector. Next, we provide a review of the butterfly factorization, but for more details we refer the
reader to \cite{BF1d,BF2d}.

\begin{figure}
\centering
 \includegraphics[width=1\textwidth,trim=1.75cm 11.65cm 1.75cm 4.25cm,clip]{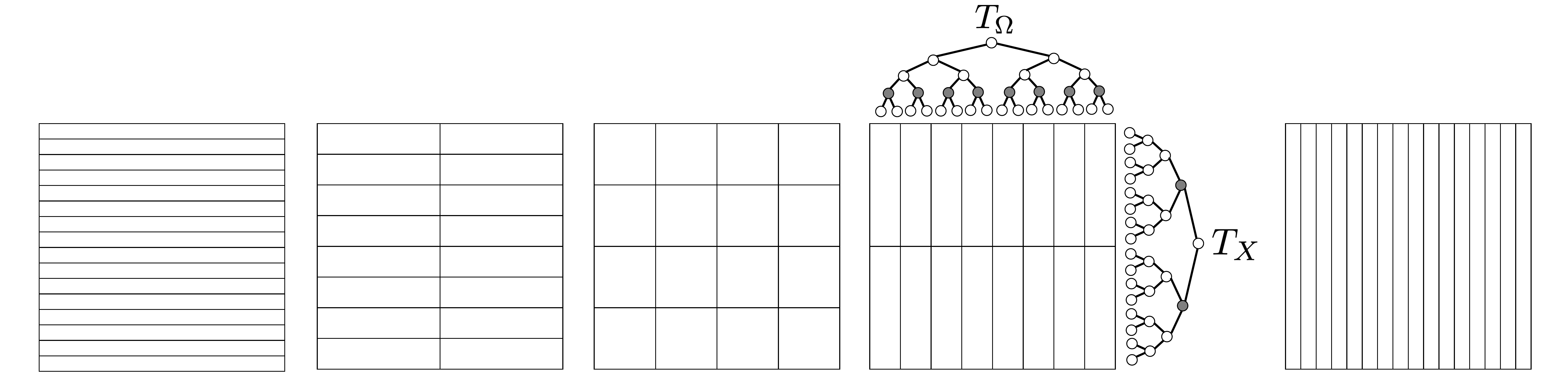}
\caption{Illustration of a matrix having complementary low-rank property, with binary partition
  trees $\Tree_X$ and $\Tree_{\Omega}$ and $L=5$. Each of the illustrations shows a different splitting
  of the matrix indices according to different levels of the trees, where each block is
  low-rank. The fourth image shows the nodes (dark nodes) used from each tree for the partition of
  the matrix, which correspond to complementary levels.}
\label{fig:1d_BFF}
\end{figure}

Define two trees $\Tree_X$ and $\Tree_{\Omega}$ of the same depth $L$, constructed by recursive
dyadic partitioning of each domain $X$ and $\Omega$. The nodes of each tree are subdomains of $X$
and $\Omega$ respectively. Level $\ell= L$ contains all the leaves of the tree and level $\ell= 0$
contains the root. With this partitioning, $\K$ satisfies the \textit{complementary low-rank
  condition} as shown in \cite{Candes}. That is for any level $0<\ell<L$, any node $\In$ at level
$\ell$ of $\Tree_X$ and any node $\Jn$ at level $L-\ell$ of $\Tree_{\Omega}$, the submatrix
$\K_{\ell;\In\Jn}$ is numerically low-rank. In this section, we use the symbol $K_{\ell;ij}$ to
represent the matrix $K$ restricted to the rows indexed by the DOFs in the $\In$-th node of
$\Tree_X$ at level $\ell$ and to the columns indexed by the DOFs in the $\Jn$-th node of
$\Tree_{\Omega}$ at level $L-\ell$.  Figure \ref{fig:1d_BFF} provides an illustration of a matrix
partitioned at different complementary levels showing different low-rank block structure due to the
complementary low-rank condition.

The butterfly factorization algorithm starts at the middle level $\ell=L/2$. Here, we assume that $L$
is even in order to ease the description of the algorithm.

\textbf{Level $\ell=L/2$.} For each pair of nodes $\In$ at level $\ell$ of $\Tree_X$ and $\Jn$ at
level $L-\ell$ of $\Tree_{\Omega}$, the algorithm computes a low-rank approximation with rank $r$ of
the submatrix $K_{\ell;ij}$,
\begin{align}\K_{\ell;\In\Jn}\approx U_{\ell;\In\Jn}D_{\ell;\In\Jn}V_{\ell;\In\Jn}^*\end{align}
with $U_{\ell;\In\Jn}\in\mathbb{C}^{N_{\In}\times r}$, $V_{\ell;\In\Jn}\in\mathbb{C}^{N_{\Jn}\times r}$, $D_{\ell;\In\Jn}\in\mathbb{C}^{r\times r}$ and $N_{\In}$ being the number of DOFs in the $i$-th node. 
One can now approximate $K$ as
\begin{align}\K\approx U_{\ell}D_{\ell}V_{\ell}^*\end{align}
where $D_{\ell}\in\mathbb{C}^{Nr\times Nr}$ is a weighted permutation matrix with only $O(N)$
non-zero entries and $U_{\ell},V_{\ell}\in\mathbb{C}^{N\times Nr}$ are block diagonal matrices. Each
diagonal block, denoted as $U_{\ell;i}$, is the horizontal concatenation of all $U_{\ell;\In\Jn}$
matrices for every node $\Jn$ at level $L-\ell$ of $\Tree_{\Omega}$. Similarly for $V_{\ell}$, each
diagonal block $V_{\ell;\Jn}$ is the horizontal concatenation of all $V_{\ell;\In\Jn}$ matrices for
every node $\In$ at level $\ell$ of $\Tree_{X}$. As a result, matrices $U_{L/2}$ and $V_{L/2}$ have
$O(N^{3/2}r)$ non-zero entries and therefore need to be further factorized to achieve a factorization that can be applied in quasi-linear
time to random vectors.

\textbf{Level $\ell=L/2,\cdots,L-1$.} At each level, $U_{\ell}$ and $V_{\ell}$ are recursively factorized by
\begin{align}
 \label{eq_BFF_nested}
    U_{\ell} = U_{\ell+1}\C_{\ell} \quad \text{and} \quad V_{\ell} = V_{\ell+1}H_{\ell}. 
\end{align}
Such factorization \eqref{eq_BFF_nested} is obtained by using nested bases, since each diagonal
block $U_{\ell;\In}$ can be obtained from the diagonal blocks $U_{\ell+1;C(\In)}$ of the children
$C(\In)$ of $\In$ in $\Tree_X$ at next level. In particular consider the 1D setting, where each node
has two children nodes except at the leaves level. One can partition $U_{\ell;\In}$ by splitting
first the matrix into two matrices indexed by DOFs of the children nodes $\In_1, \In_2$ of node $\In$, and
then merging together the subblocks of the form $U_{\ell;\In_k\Jn}$ that share the same parent on
$\Tree_{\Omega}$. This is, for each node $i$, one can define blocks of the form $U_{\ell;\hat{\imath}\hat{\jmath}} =
[\U_{\ell;\hat{\imath}\Jn},U_{\ell;\hat{\imath}(\Jn+1)}]$, where $\Jn$ and $\Jn+1$ are nodes in
level $L-\ell$ of $\Tree_{\Omega}$ with common parent node $\hat{\jmath}$, and $\hat{\imath}$ is
a child of node $\In$. This partitioning of the form
\begin{align}
    U_{\ell;\In} = \left[{\begin{array}{c|cc|cc|c}
    \cdots& U_{\ell;\In_1\Jn} & U_{\ell;\In_1(\Jn+1)} & U_{\ell;\In_1(\Jn+2)} & U_{\ell;\In_1(\Jn+3)} &\cdots \\
    \hline
    \cdots& U_{\ell;\In_2\Jn} & U_{\ell;\In_2(\Jn+1)} & U_{\ell;\In_2(\Jn+2)} & U_{\ell;\In_2(\Jn+3)} &\cdots \\
    \end{array}}\right],
\end{align}
can analogously be extended to the 2D setting, with quadtrees $\Tree_X$ and $\Tree_{\Omega}$, by
partitioning the matrix $U_{\ell;\In}$ over 4 children instead of 2 from $\Tree_X$, and
concatenating each 4 sibling nodes from $\Tree_\Omega$.  Notice that due to the complementary
low-rank condition, each $U_{\ell;\hat{\imath}\hat{\jmath}}$ is numerically low-rank and therefore
one can obtain the corresponding low-rank approximation
\begin{align}
 \label{eq_BFF_child_parent}
U_{\ell;\hat{\imath}\hat{\jmath}} \equiv 
[\U_{\ell;\hat{\imath}\Jn},U_{\ell;\hat{\imath}(\Jn+1)}] =  U_{\ell+1;\hat{\imath}\hat{\jmath}} \C_{\ell;\hat{\imath}\hat{\jmath}},
\end{align}
where $U_{\ell+1;\hat{\imath}\hat{\jmath}}\in\mathbb{C}^{N_{\hat{\imath}}\times r'}$,  $U_{\ell;\hat{\imath}\hat{\jmath}}\in\mathbb{C}^{N_{\hat{\imath}}\times 2r}$ with $r'<2r$ and $\C_{\ell,\hat{\imath}\hat{\jmath}}$ is a matrix of size $r' \times 2r$ for 1D FIOs, and $r'
\times 4r$ for 2D FIOs. Rearranging all these
factors for all nodes $\hat{\imath}$ and $\hat{\jmath}$ in a block-diagonal matrix
$U_{\ell+1}\in\mathbb{C}^{N\times Nr'}$ and a sparse matrix $\C_{\ell}$ with only $O(N)$ non-zero entries,
one obtains the factorization from \eqref{eq_BFF_nested}. The same recursion can be applied to
$V_{\ell}$.

By applying \eqref{eq_BFF_nested} recursively, one obtains at level $\ell=L$ the approximate
\textit{butterfly factorization} $\Ka$ of $K$ as
\begin{align}
  \label{eq_BFF_fact}
  \K \approx \Ka \equiv U_L\C_{L-1} \cdots \C_{L/2}D_{L/2}H_{L/2}^* \cdots H_{L-1}^*V_{L}^*
\end{align}
with $O(N^{1.5})$ factorization cost. Since all the factors in \eqref{eq_BFF_fact} have $O(N)$
non-zero entries, the butterfly factorization of $K$ can be applied in quasi-linear time $O(N\log
N)$ to a random vector. Additionally, one can easily approximate the transpose of $\K$
\[
K^* \approx \Ka^* \equiv V_L H_{L-1} \cdots H_{L/2}D_{L/2}^*\C_{L/2}^* \cdots \C_{L-1}^*U_{L}^*.
\]
Notice that most of the matrices in \eqref{eq_BFF_fact} are rectangular and therefore the
factorization is not readily invertible.

\section{Hierarchical matrix construction}\label{sec:Peeling}

Since $\K$ is a square matrix, its inverse is equal to its pseudo-inverse $\K^{-1}=(\Kt\K)^{-1}\Kt$. In this section, in order to build $K^+$, we are interested in a representation of the matrix $\Kt\K$ that can be constructed and inverted in quasi-linear time. In
particular, we propose to use the $\mathcal{H}$-matrix framework with different admissibility conditions
depending on the dimension of the problem.

\subsection{Hierarchical matrix approximation of $\Kt\K$}

In order to analyze the properties of $\Kt\K$ we start by decomposing the FIO from \eqref{eq1} as
$\mathcal{L} = \mathcal{K} \F$, where $\F(f) = \hat{f}$ is the Fourier transform, and
\begin{align}
  \label{eq_S}
  (\mathcal{K} f)(x) = \int_{\Rb^d} a(x,\xi)e^{2\pi i \Phi(x,\xi)}f(\xi)\text{d}\xi.
\end{align}

Defining the adjoint operator of $\mathcal{K}$ as
\begin{align}
  \label{eq_adj}
  (\mathcal{K}^* f)(\xi) = \int_{\Rb^d} \bar{a}(x,\xi)e^{-2\pi i \Phi(x,\xi)}f(x)\text{d}x,
\end{align}
with $\bar{a}(x,\xi)$ the complex conjugate of $a(x,\xi)$, we can construct the operator $\mathcal{K}^*\mathcal{K} $ defined by
\begin{align}
  \label{eq_normal}
  (\mathcal{K}^*\mathcal{K} f)(\xi) = \int_{\mathbb{R}^d} \adjkernel(\xi,\eta)f(\eta)\text{d}\eta,
\end{align}
where $\adjkernel$ is the kernel of the operator $\mathcal{K}^*\mathcal{K}$ defined by
\begin{align}
  \label{eq_kernel}
  \adjkernel(\xi,\eta) = \int_{\Rb^d} \bar{a}(x,\xi)a(x,\eta)e^{2\pi i\left[ \Phi(x,\eta)-\Phi(x,\xi)\right]}\text{d}x.
\end{align}
The following proposition gives bounds for the kernel $\adjkernel$ (see Chapter IX Section 3.1.2 of
\cite{Stein}).

\begin{proposition}\label{prop:hodlr}
Under the assumptions that\begin{itemize}
    \item $a(x,\xi)$ is a fixed smooth function of compact support in $x$ and $\xi$,
    \item the phase function $\Phi(x,\xi)$ is real-valued, homogeneous of degree 1 in $\xi$ and smooth in $(x,\xi)$ for $\xi\neq0$ on the support of $a$,
    \item $\Phi$ satisfies the non-degeneracy condition for $\xi\neq0$:
det$\Big(\frac{\partial^2\Phi}{\partial x_i\partial\xi_i}\Big)\neq0$,
\end{itemize}
we have that
\[|\adjkernel(\xi,\eta)|\leq A_k(1+|\xi-\eta|)^{-k}, \quad \forall k\geq0 \]
with $A_k$ a constant depending on $k$.
\end{proposition}

The discrete analog of $\mathcal{K}^*\mathcal{K}$ for the spatial and frequency domains
$\mathcal{X}$ and $\Omega$ on the Cartesian grids \eqref{eq:X_domain} and \eqref{eq:freq_domain} is
\begin{align} \label{eq_normal_discr}
  S \equiv K^*K
\end{align}
where $K\in\mathbb{C}^{N\times N}$ is the discrete FIO matrix from \eqref{eq6} and $S$ has been overwritten to represent the discrete matrix for the kernel $S(\xi,\eta)$ of $\mathcal{K}^*\mathcal{K}$.
According to Proposition \ref{prop:hodlr}, in the continuous case the kernel $\adjkernel(\xi,\eta)$
decays fast while in the discrete case we observe experimentally that the entries of the hermitian
matrix $S\equiv\K^*\K$ also decay fast, so that $S\equiv \K^*\K$ can be represented as a banded
matrix. Additionally, we observe that the interaction between non-adjacent clusters of points is
numerically low-rank. Therefore, one can represent $S\equiv\K^*\K$ using hierarchical matrices. In
particular, we use HODLR \cite{Aminfar2} for the 1D setting and $\mathcal{H}$-matrices
\cite{Hackbusch} format for the 2D setting, using weak and strong admissibility respectively.

The admissibility condition is defined as follows: two disjoint sets of indices $\I$ and $\J$ are admissible if
\begin{align}
  \min(\text{diam}(\I), \text{diam}(\J)) \leq \alpha \cdot \text{dist}(\I , \J)
\end{align}
for $\alpha>0$, where $\text{diam}(\I) = \max_{i,j\in \I} |x_i -
x_j|$ and $\text{dist}(\I,\J) = \min_{i\in \I,j\in \J} |x_i - x_j|$. For well-separated
(admissible) clusters of points, the corresponding submatrix block in $K^*K$ is low-rank.

If weak admissibility is used, $\Kt\K$ can be represented with HODLR format, i.e. every off-diagonal
block of the matrix is considered low-rank. Otherwise, strong admissibility is used with $\alpha=1$,
leading to the $\mathcal{H}$-matrix format, where the interactions with the adjacent clusters are considered
full-rank while the interactions with the remaining clusters are low-rank. For any subblock of the
matrix, we consider it is numerical low-rank if its numerical rank increases very slowly with
$N$, no more quickly than $O(\log N)$.

The pseudo-inverse of $K$ can be computed as $\K^{+} = (\Kt\K)^{-1}\Kt$ or as $\K^{+} =
\Kt(\K\Kt)^{-1}$. However, for the examples of FIO from Section \ref{sec:NR}, we
observe experimentally that while both $\Kt\K$ and $\K\Kt$ can be approximated with hierarchical matrices, the ranks
of the admissible blocks are lower for the former. Therefore, we propose to construct the inverse of
$\K$ following $\K^{-1} = (\Kt\K)^{-1}\Kt$.  The observed experimental asymptotic complexity of the
ranks of the admissible blocks for the examples in Section \ref{sec:NR} are reported in Table
\ref{table:ranks}. This experimental observations justify the use of weak admissibility for 1D FIOs to
represent $\Kt\K$ as a hierarchical matrix in HODLR format. However, we need to use strong
admissibility and the $\mathcal{H}$-matrix format for 2D FIOs to avoid ranks scaling with
$O(N^{1/2})$.

\begin{table}[htb]
\centering

\begin{tabular}{cc|c|c}
\toprule
&&\multicolumn{2}{c}{Admissibility} \\
\midrule
&&Weak & Strong\\
\midrule
\multirow{2}{*}{1D} & $\K\Kt$ & $O(\log N)$ & $O(1)$\\
& $\Kt\K$ & $O(\log N)$ & $O(1)$\\
\midrule
\multirow{2}{*}{2D} & $\K\Kt$ & $O(N^{2/3})$ & $O(N^{1/2})$\\
& $\Kt\K$ & $O(N^{1/2})$ & $O(\log N)$\\
\bottomrule
\end{tabular}
\caption{Approximate observed rank scalings of admissible blocks from the composition of a discrete FIO matrix with its adjoint for FIOs used in section \ref{sec:NR}, where strong admissibility is considered with $\alpha = 1$.}
\label{table:ranks}
\end{table}

\subsection{Review of peeling algorithm for 1D HODLR matrices}

In order to construct the HODLR representation of $S\equiv \Kt\K$ we use the peeling algorithm
\cite{Peeling}, which only requires matrix-vector multiplications that can be performed rapidly in
quasi-linear time with the butterfly factorization $\Ka\approx\K$.


For a given matrix $A$, we start by defining a uniform cluster tree $\Tree_{\Omega}$ with $L$ levels, that recursively
partitions the domain $\Omega$ into $2^{\ell}$ subdomains at each level
$\ell=0,1,...,L$, with level $\ell=0$ being the root level and level $\ell=L$ containing the leaves
of the tree. To obtain the HODLR representation of $A$, at every level each subdomain is split into two child subdomains, and its
interaction block in $A$ is approximated by a low-rank matrix.

Following the peeling algorithm from \cite{Peeling}, we can obtain the low-rank approximations of
submatrices at each level $\ell$, starting from level $\ell=1$ up to the leaves level $\ell = L$.
Here, we consider obtaining the HODLR approximation for an hermitian matrix $\A \equiv \Ka^* \Ka$.  For
the rest of the section, we denote $\A_{\ell;ij}$ the subblock of the matrix $\A$ restricted on the
rows and columns indexed by the $i$-th and $j$-th nodes at level $\ell$ of the tree $\Tree_{\Omega}$
respectively.

\textbf{Level $\ell=1$}. At this level, the domain is divided into two subdomains and the matrix takes the form
\begin{align}
\A = \begin{bmatrix}
  \A_{1;11} &\A_{1;12} \\
  \A_{1;21}&\A_{1;22} \\
\end{bmatrix},
\label{eqn:matrix_hodlr}
\end{align}
where $\A_{1;12}$ and $\A_{1;21}$ can be approximated by a low-rank matrix with rank $k$.

One starts by generating two random matrices $R_{1;1}$ and $R_{1;2}$ of size $N/2\times(k+c)$, with $c = O(1)$ an oversampling parameter. One can perform the following matrix-matrix multiplications to obtain
\begin{align}
   \begin{bmatrix}
  \A_{1;11} &\A_{1;12} \\
  \A_{1;21} &\A_{1;22} \\
\end{bmatrix}\begin{bmatrix}
  R_{1;1} \\
  0 \\
\end{bmatrix} = \begin{bmatrix}
 \A_{1;11}R_{1;1}\\
  \A_{1;21}R_{1;1} \\
\end{bmatrix},\\
 \begin{bmatrix}
  \A_{1;11} &\A_{1;12} \\
  \A_{1;21}&\A_{1;22} \\
\end{bmatrix}\begin{bmatrix}
 0 \\
  R_{1;2} \\
\end{bmatrix} = \begin{bmatrix}
\A_{1;12}R_{1;2}\\
  \A_{1;22}R_{1;2} \\
\end{bmatrix}.
\label{eqn:peeling_lvl1_2}
\end{align}

Since $A$ is hermitian, $\A_{1;12} = \A_{1;21}^*$. Now, using QR factorization (or alternatively SVD), we can obtain a unitary column base matrix $U_1$ of $\A_{1;21}R_{1;1}$ and a unitary column base matrix $U_2$ of $\A_{1;21}^*R_{1;2}$. A low-rank approximation \cite{Halko} of $\A_{1;21}$ is obtained as
\begin{align}
 \A_{1;21} \approx \hat{\A}_{1;21}= U_1 (R_{1;2}^*U_1)^{+}(R_{1;2}^*\A_{1;21}R_{1;1})(U_2^*R_{1;1})^{+}U_2^*.
\label{eqn:lvl1_LR}
\end{align}
One can now construct the matrix $\hat{\A}_{1}$ containing the low-rank approximations at this level as
\begin{align}
 \hat{\A}_{1} =\begin{bmatrix}
 0 & \hat{\A}_{1;21}^*\\
  \hat{\A}_{1;21} & 0\\
\end{bmatrix}.
\label{eqn:lvl1_LRF}
\end{align}

\textbf{Level $\ell=2$}. 
Generate the matrix $[R_{2;1}^T, 0, R_{2;3}^T, 0]^T$ where $R_{2;1}$ and $R_{2;3}$ are random matrices of size $N/4\times(k+c)$. In order to obtain the low-rank approximations of $\A_{2;21}$ and $\A_{2;43}$ at level $\ell=2$, we can now multiply $A$ with the random matrices and subtract the contribution of the already computed low-rank approximations at previous levels,
\begin{align}
\begin{split}
A [R_{2;1}^T, 0, R_{2;3}^T, 0]^T - \hat{\A}_{1}[R_{2;1}^T 0 R_{2;3}^T 0]^T&=\\ \begin{bmatrix}
  \A_{2;11} &\A_{2;12}& \multicolumn{2}{c}{\multirow{2}{*}{$\A_{1;12}-\hat{\A}_{1;12}$}}\\
  \A_{2;21} &\A_{2;22} &&\\
 \multicolumn{2}{c}{\multirow{2}{*}{$\A_{1;21}-\hat{\A}_{1;21}$}}&\A_{2;33} &\A_{2;34}\\
  &&\A_{2;43} &\A_{2;44} \\
\end{bmatrix}\begin{bmatrix}
  R_{2;1} \\
  0 \\
   R_{2;3} \\
  0 \\
\end{bmatrix} &=\\ 
\begin{bmatrix}
 \begin{pmatrix}
\A_{2;11}R_{2;1}\\
\A_{2;21}R_{2;1}
\end{pmatrix} + (\A_{1;12}-\hat{\A}_{1;12})\begin{pmatrix}
R_{2;3}\\
0
\end{pmatrix}\\
  \begin{pmatrix}
\A_{2;33}R_{2;3}\\
\A_{2;43}R_{2;3}
\end{pmatrix} + (\A_{1;21}-\hat{\A}_{1;21})\begin{pmatrix}
R_{2;1}\\
0
\end{pmatrix}\\
\end{bmatrix} &\approx 
\begin{bmatrix}
\A_{2;11}R_{2;1}\\
\A_{2;21}R_{2;1}\\
\A_{2;33}R_{2;3}\\
\A_{2;43}R_{2;3}
\end{bmatrix}.
\end{split}
\label{eqn:peeling_lvl2}
\end{align}
Similarly, one can perform the same operation with $(0, R_{2;2}^T, 0, R_{2;4}^T)^T$ to obtain $\A_{2;12}R_{2;2}$ and $\A_{2;34}R_{2;4}$. The low-rank approximations for $\A_{2;21}$ and $\A_{2;43}$ are computed as in \eqref{eqn:lvl1_LR}, resulting
in the low-rank approximation matrix at this level
\begin{align}
 \hat{\A}_{2} =\begin{bmatrix}
 0 & \hat{\A}_{2;21}^*\\
  \hat{\A}_{2;21} & 0\\
   &&0 & \hat{\A}_{2;43}^*\\
  &&\hat{\A}_{2;43} & 0\\
\end{bmatrix}.
\label{eqn:lvl2_LRF}
\end{align}

\textbf{Levels $\ell=3,\cdots,L$}.
Generate two matrices $R_1,R_2\in\mathbb{C}^{N\times (k+c)}$  containing the random matrices for each of the $2^\ell$ numbered subdomains as 
\begin{align}R_1(\I_i,:) = \begin{cases} 0, & i\text{ is even}\\  
    R_{\ell;i},  & \text{otherwise}
    \end{cases}, \qquad R_2(\I_i,:) = \begin{cases} R_{\ell;i} & i\text{ is even}\\  
    0,  & \text{otherwise}
    \end{cases}
    \end{align}
with $\I_i$ being the DOFs of the $i$-th node of $\Tree_\Omega$ at level $\ell$ for $i=1,\cdots,2^{\ell}$.

Similar to level $\ell=2$, we compute $\A R_j-\sum_{l=1}^{\ell-1}\hat{\A}_lR_j$ with $j=1,2$, to obtain $\A_{\ell;(i+1)i}R_{\ell;i}$ and $\A_{\ell;(i+1)i}^*R_{\ell;i+1}$ for each node $i$ such that $(i\mod2)=1$. 
Analogously to the two previous levels, one can compute the low-rank approximations of each non-zero subblock $\hat{\A}_{\ell;(i+1)i} \approx \A_{\ell;(i+1)i}$ of $\hat{\A}_{\ell}$.

At the last level $\ell=L$, one can compute the block diagonal matrix $D$ which approximate the diagonal
full-rank blocks of matrix $\A$ for each of the $2^{L}$ leave nodes of
$\Tree_{\Omega}$. The HODLR matrix representation of $\A$ can be expressed as
\begin{align}
   \Sa \equiv D + \sum_{\ell=1}^{L}\hat{A_{\ell}}.
\end{align}

Since $A\equiv \Ka^*\Ka \approx K^*K \equiv S$, $\Sa$ is an approximation to $S$.

\subsection{Review of peeling algorithm for 2D $\mathcal{H}$-matrix}

For 2D FIO matrices we use the $\mathcal{H}$-matrix format to approximate $\Ka^*\Ka$, which only uses
low-rank approximations on the interactions between non-intersecting and non-adjacent clusters. In
this section, we review the algorithm used to construct such representation.

We start by defining a uniform quad-tree $\Tree_{\Omega}$ with $L$ levels, that partitions the domain $\Omega$ into $2^{\ell}\times2^{\ell}$ square  cells  at  each  level $\ell=0,1,...,L$. Therefore, at every level each subdomain is split into four child subdomains. The leaves  correspond  to  level $\ell=L$ and  the  root  to  level $\ell=0$. Each cell at level $\ell$ is represented as $I_{\ell;ij}$ with $i=1,\cdots,2^{\ell}$ and $j=1,\cdots,2^{\ell}$.
For each node/cell $I$ of the tree we define the following two lists:
\begin{itemize}
    \item \textit{NL($I$): neighbor list of cell $I$}. This is the list containing the cell $I$ and the adjacent cells on level $\ell$. 
    \item \textit{IL($I$): interaction list of cell $I$}. This is the list containing all the cells on level $\ell$ that are children of cells in NL($P(I)$) minus the cells in NL($I$), where $P(I)$ is the parent cell of $I$.
\end{itemize}

The peeling algorithm from \cite{Peeling} used to obtain the $\mathcal{H}$-matrix representation for the 2D case is similar to the 1D case, starting from the root level and going all the way to the leaves, with the following main differences:
\begin{itemize}
    \item At every level $\ell$, for each cell $I$ the algorithm only computes the low-rank approximation of subblocks $\A_{\ell;\I\J}$ of $A$ restricted to the rows indexed by the set of indices $\I$ of DOFs in cell $I$ and the columns indexed by the set of indices $\J$ of DOFs in cell $J$, with  $J\in IL(I)$.
    
    \item In order to compute the low-rank approximation of $\A_{\ell;\I\J}$ one introduces the following 64 sets of cells $\mathcal{P}_{pq}$ for $1\leq p,q\leq 8$ with
    \begin{align}
    \mathcal{P}_{pq} = \{I_{\ell;ij} | i\equiv p \bmod 8, j\equiv q \bmod 8\}.
    \label{eq_spq}
    \end{align}
    One can now construct 64 proving matrices, one for each set in $\mathcal{P}_{pq}$, of the form
    \[R(\I,:) = \begin{cases} R_I, & I\in \mathcal{P}_{pq}\\  
    0,  & \text{otherwise,}
    \end{cases}\]
    with $R_I$ a random matrix of size $|\I|\times(k+c)$.
    This provides for every cell $I$ and cell $J\in$ NL($I$) the matrices $\A_{\ell;\I\J}R_J$ and $\A_{\ell;\J\I}R_I$, which allows to efficiently construct a low-rank approximation
    \[\A_{\ell;\I\J}\approx U_{\ell;\I\J}\B_{\ell;\I\J}U_{\ell;\J\I}^*,\] 
    with $U_{\ell;\I\J}\in\mathbb{C}^{|\I|\times k}$ and $\B_{\ell;\I\J}\in\mathbb{C}^{k\times k}$. 
\end{itemize}

\begin{figure}
\centering
 \includegraphics[width=0.25\textwidth,trim=1.65cm 1.65cm 0.75cm 1.15cm,clip]{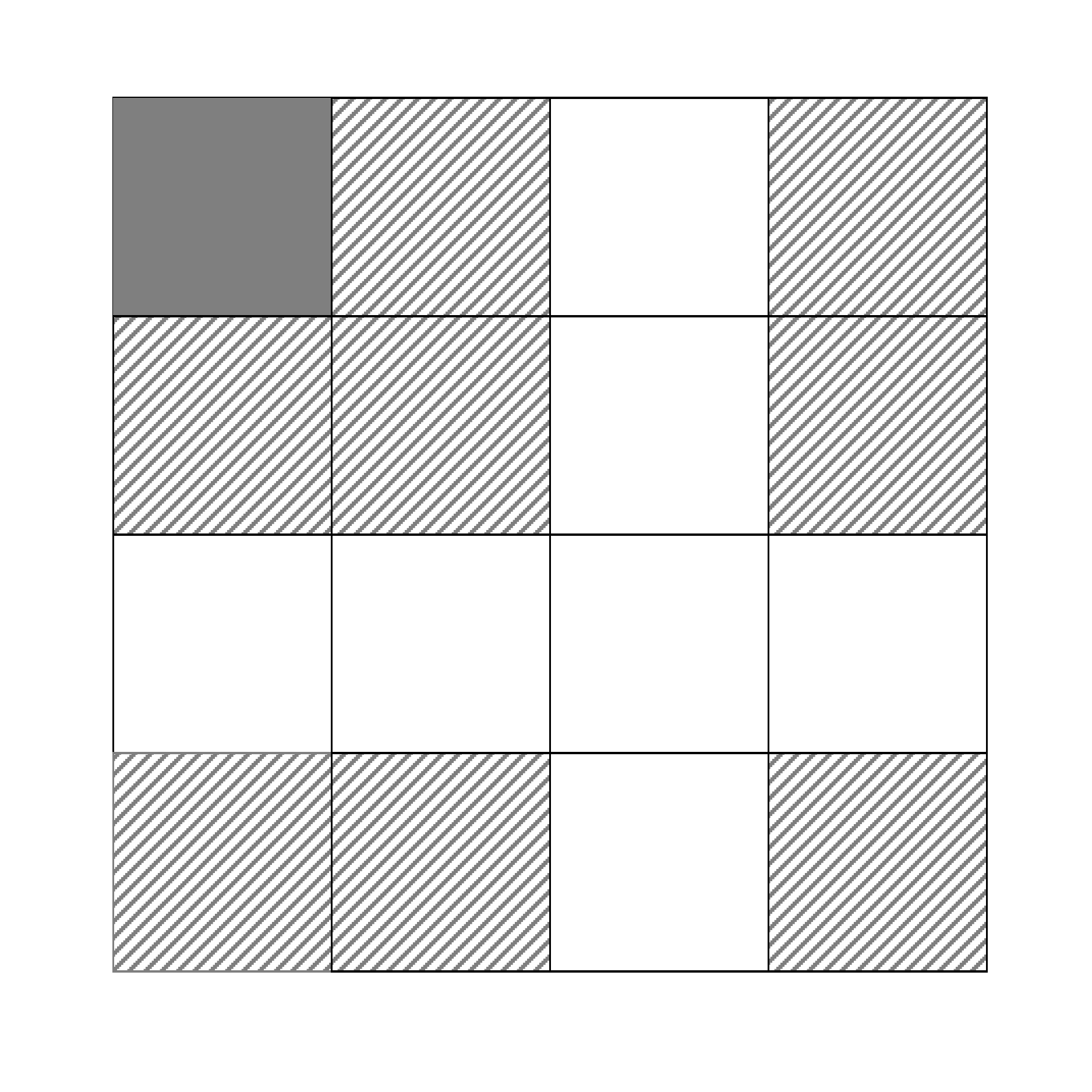}
\caption{Square cells at level $\ell=2$ of domain $\Omega$. A particular cell $I$ (with solid grey filling) at level $\ell=2$ and the cells in $NL(I)$ (with dashed grey filling).}
\label{fig:2d_Hmatrix}
\end{figure}

For more details we refer the reader to \cite{Peeling}. If we use $\A\equiv \Ka^*\Ka \approx \K^*\K \equiv S$, the matrix $\Sa$ obtained from the peeling algorithm applied to $\A$ is an approximation of $S$. A particularity of constructing the
$\mathcal{H}$-matrix for the 2D FIO examples in Section \ref{sec:NR} is that there is periodicity in
both directions in the domain $\Omega$ when accounting for the admissible blocks. What that means is
that any subblock of $\Ka^*\Ka$ restricted to the rows of cell $I$ and the columns indexed by any
other cell $J$ that becomes adjacent to $I$ after padding the domain $\Omega$ with itself on both
directions, is not considered low-rank. Therefore, $J$ is added to NL$(I)$ in order to avoid
numerical ranks increasing quickly with $N$. An example of a particular cell $I$ and its neighbor
list NL($I$) is illustrated in Figure \ref{fig:2d_Hmatrix}.

\section{Hierarchical matrix inversion}\label{sec:HIF}
In this section we propose an algorithm to invert matrices expressed in hierarchical format using
recursive skeletonization (RSS) \cite{Ho} and the hierarchical interpolative decomposition (HIF)
\cite{HIFIE}.

\subsection{1D case}
In the 1D setting, we have constructed an HODLR approximation $\Sa \approx S\equiv K^*K$. One can
easily invert $\Sa$ using recursive skeletonization \cite{Ho}, leading to an algorithm that resembles the
work in \cite{XiaHSS10}.

In this subsection, we use $A$ to denote the HODLR approximation $\Sa$, which takes the form 
\begin{align}
  \A =\begin{bmatrix}
  \Hm_{2;11} & \U_{2;1}\B_{2;21}^*\U_{2;2}^*& \multicolumn{2}{c}{\multirow{2}{*}{$\U_{1;1}\B_{1;21}^*\U_{1;2}^*$}}\\
  \U_{2;2}\B_{2;21}\U_{2;1}^* &\Hm_{2;22} &&\\
  \multicolumn{2}{c}{\multirow{2}{*}{$\U_{1;2}\B_{1;21}\U_{1;1}^*$}}&\Hm_{2;33} &\U_{2;3}\B_{2;43}^*\U_{2;4}^*\\
  &&\U_{2;4}\B_{2;43}\U_{2;3}^* &\Hm_{2;44}\\
  \end{bmatrix}
\end{align}
with $U_{\ell;i}\in \mathbb{C}^{N/2^{\ell}\times r}$ being the low-rank orthogonal column base of
$\Hm_{\ell;i(i+1)}$, the subblock of matrix $\Hm$ restricted by the DOFs of cells $i$ and $i+1$ at
level $\ell$, and $B_{\ell;i(i+1)}\in \mathbb{C}^{r\times r}$. The diagonal blocks $\Hm_{\ell;ii}$
are recursively factorized following the levels of the tree up to the leaves, by making the
off-diagonal matrix terms low-rank.

The algorithm to invert $A$ keeps eliminating nodes at each level starting from the leaves up to the
root. Let $A_{\ell}$ be the state of the matrix $A$ at level $\ell$, with $A_L=A$, and let
$\S_{\ell}$ be the \textit{active} DOFs at level $\ell$, i.e. DOFs that have not been
decoupled/eliminated in previous levels throughout the factorization. To ease the description of the
algorithm, we will drop $\ell$ when representing subblocks of matrix $A_{\ell}$ restricted to rows
indexed by $\I$ and columns indexed by $\J$, $A_{\I\J}\equiv(A_{\ell})_{\I\J}$. The inversion
algorithm proceeds from the leaves level $\ell=L$ of the tree $\Tree_{\Omega}$ to level $\ell=1$
performing the following two steps at each level $\ell$:

\begin{itemize}
\item \textbf{Skeletonization of nodes at level $\ell$.}
  There are $q_{\ell} = 2^\ell$ nodes at level $\ell$. For each node we let $\I$ be the set of indices for the \textit{active} DOFs associated with that node. First, we construct the matrix 
  \begin{align}
    \bar{U}_{\ell;\I} =     \begin{bmatrix}
      U_{\ell;\I}^*\\
      \vdots\\
      U_{1;\I}^*
    \end{bmatrix},
  \end{align}
  where $U_{j;\I}$ are the corresponding columns of the low-rank base from the HODLR representation at level $j$ of 
  $\Tree_{\Omega}$ for indices $\I$ in the corresponding node. Assume that $\bar{U}_{\ell;\I}$ can be approximated by a low-rank matrix with relative precision $\epsilon$ and rank $k$. Using the interpolative decomposition (ID) \cite{ID}, we can find a disjoint
  partitioning $\I=\Ir\cup\Is$ into \textit{redundant or fine} ($\Ir$) and \textit{skeleton or
    coarse} ($\Is$) DOFs such that
  \begin{align}\label{eq_id}
\bar{U}_{\ell;\Ir} = \bar{U}_{\ell;\Is}\T_{\I} + E_{\I}, \quad \|E_{\I}\| = O(\epsilon\|\bar{U}_{\ell;\I}\|),
\end{align}
where $\T_{\I}\in\mathbb{R}^{k\times(|\I|-k)}$ is an \textit{interpolation matrix}, which serves the purpose of approximating the redundant columns $\bar{U}_{\ell;\Ir}$ by a linear combination of the skeleton columns $\bar{U}_{\ell,\Is} $. Additionally, the interpolation matrix is constructed using strong rank-revealing QR so that $\|\T_{\I}\|$ is not too large.
Using this decomposition we can approximately zero out the redundant DOFs in the low-rank bases of the form $U_{\ell;\I}^*$ resulting, up to a permutation, in
\begin{align}
    \bar{U}_{\ell;\I}
    \begin{bmatrix}
           I & \\
           -\T_{\I} & I 
    \end{bmatrix} \approx \begin{bmatrix}
         0&U_{\ell;\Is}^*\\
           \vdots&\vdots\\
          0&U_{1;\Is}^*
    \end{bmatrix}.
    \label{eq:ID_1d}
    \end{align}
    
    Notice that zeroing out the columns $\Ir$ in \eqref{eq:ID_1d}, also approximately zeroes out the redundant columns $\Ir$ of the subblock $\A_{\R\I}$ of the matrix $\A_{\ell}$, with $\R=\S_{\ell}\setminus\I$.
    By symmetry of $\A$, the rows $\Ir$ of the subblock $\A_{\I\R}$ are also zeroed out.
    Introducing the \textit{zeroing matrix} $\Z_{\I}$ and applying it on both sides leads, up to a permutation, to
    \begin{align}
    \label{eq_sk}
      \Hs_{\ell} = \Z_{\I}^*\A_{\ell} \Z_{\I}\approx \left[{\begin{array}{cc|c}
        \Hs_{\Ir\Ir}&\Hs_{\Is\Ir}^*& \vspace*{0.1cm}\\\vspace*{0.1cm} \Hs_{\Is\Ir}&\Hm_{\Is\Is}&\Hm_{\R\Is}^*\\ \hline
        &\Hm_{\R\Is}&\Hm_{\R\R}\\
    \end{array}}\right],
  \quad
  \Z_{\I}=\left[{\begin{array}{cc|c}
        I&&\\
        -\T_{\I}&I&\\ \hline
        &&I\\
    \end{array}}\right],
\end{align}
where only the terms $\AE_{\Ir\Ir}$ and $\AE_{\Is\Ir}$ are updated with
\begin{align}
\Hs_{\Ir\Ir}=\Hm_{\Ir\Ir}-\T_{\I}^{*}\Hm_{\Is\Ir}-\Hm_{\Is\Ir}^*\T_{\I}+\T_{\I}^*\Hm_{\Is\Is}\T_{\I},
\quad
\Hs_{\Is\Ir}=\Hm_{\Is\Ir}-\Hm_{\Is\Is}\T_{\I}.
\end{align}
We can now decouple the redundant DOFs $\Ir$ by performing standard Gaussian elimination
\begin{align} \label{eq_sk_3}
\E_{\Ir}^{*}\Hs_{\ell} \E_{\Ir}=\begin{bmatrix}
  I&&\\
  &\He_{\Is\Is}&\Hm_{\R\Is}^*\\
  &\Hm_{\R\Is}&\Hm_{\R\R}\\
\end{bmatrix}, \quad\E_{\Ir}=
  \begin{bmatrix}
    (\L_{\Ir}^*)^{-1} & -\Hs_{\Ir\Ir}^{-1} \Hs_{\Is\Ir}^*\\
    &I\\
    &&I
  \end{bmatrix} 
\end{align}
where $\He_{\Is\Is}$ has been updated to
$\He_{\Is\Is}=\Hm_{\Is\Is}-\Hs_{\Is\Ir}\Hs_{\Ir\Ir}^{-1}\Hs_{\Is\Ir}^*$.  At this point we have
decoupled the redundant DOFs $\Ir$ from the rest of DOFs. This process is called
\textit{skeletonization} of $\I$.  Notice that in order to build the 
factorization in quasi-linear time, it is essential to perform ID on $\bar{U}_{\ell,\I}$ instead of the block
$\A_{\R\I}$, since it only has $O(r \log N)$ rows, with $r$ the average rank of
the off-diagonal blocks at each level.

It is important to remark that performing skeletonization on $\I$ decouples its redundant DOFs
(i.e. zeroes out $\Hm_{\Ir\R}$ and $\Hm_{\R\Ir}$), while not changing the interactions of its
skeleton DOFs with the remaining DOFs, since only the block $\Hm_{\Is,\Is}$ is changed. Therefore,
the low-rank approximation of the remaining active DOFs at next levels remains the same and so does
$U_{\ell;\Is}$ for all successive levels. This might be useful if one wants to select skeletons a
priori based on their geometric interpretation in the matrix which is not lost by maintaining the
original approximate values of the off-diagonal blocks in the matrix.
     
Let $\{ \I_{\ell, i} \}_{i = 1}^{q_\ell}$ be the collection of disjoint index sets corresponding to
each node $i$ at level $\ell$. Skeletonization on all the index sets then gives
\begin{align} \label{eq_sk_4}
  \Hm_{\ell-1}\approx R_{\ell}^{*}\Hm_{\ell} R_{\ell}, \qquad
  R_{\ell}=\prod_{i=1}^{q_\ell}\Z_{\I_{\ell,i}}\E_{\Ir_{\ell,i}}.
\end{align}
The remaining active DOFs in $\Hm_{\ell-1}$ are $\S_{\ell-1} = \cup_{i=1}^{q_{\ell}}\Is_{\ell,i}$. The DOFs indexed by
$\cup_{i=1}^{q_{\ell}}\Ir_{\ell,i}$ have been decoupled, i.e. their corresponding matrix block in $\Hm_{\ell-1}$
is the identity matrix.
\item \textbf{Merge child blocks.}  After the skeletonization step, we merge the active DOFs of
  sibling nodes in the tree and move to the next level of the tree, i.e. to the parent node at level
  $\ell-1$.
\end{itemize}

\begin{figure}
\centering \captionsetup[subfigure]{labelformat=empty}
\begin{subfigure}{0.27\textwidth}
 \centering
 \includegraphics[width=\textwidth,trim=4.25cm 8.25cm 4.2cm 8.5cm,clip]{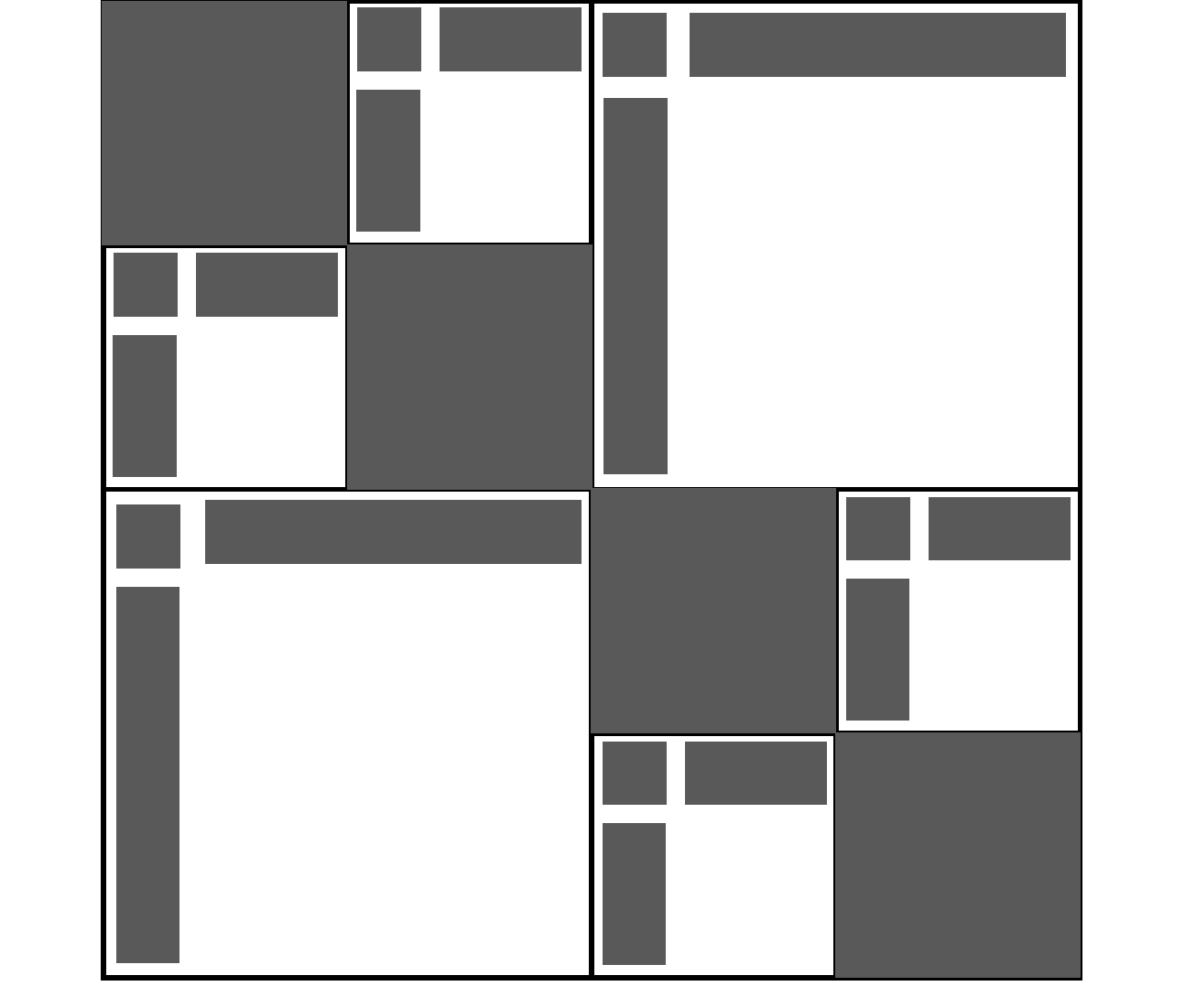}
 \caption{$\A$}
\end{subfigure}
\begin{subfigure}{0.27\textwidth}
 \centering
 \includegraphics[width=\textwidth,trim=4.25cm 8.4cm 4.2cm 8.1cm,clip]{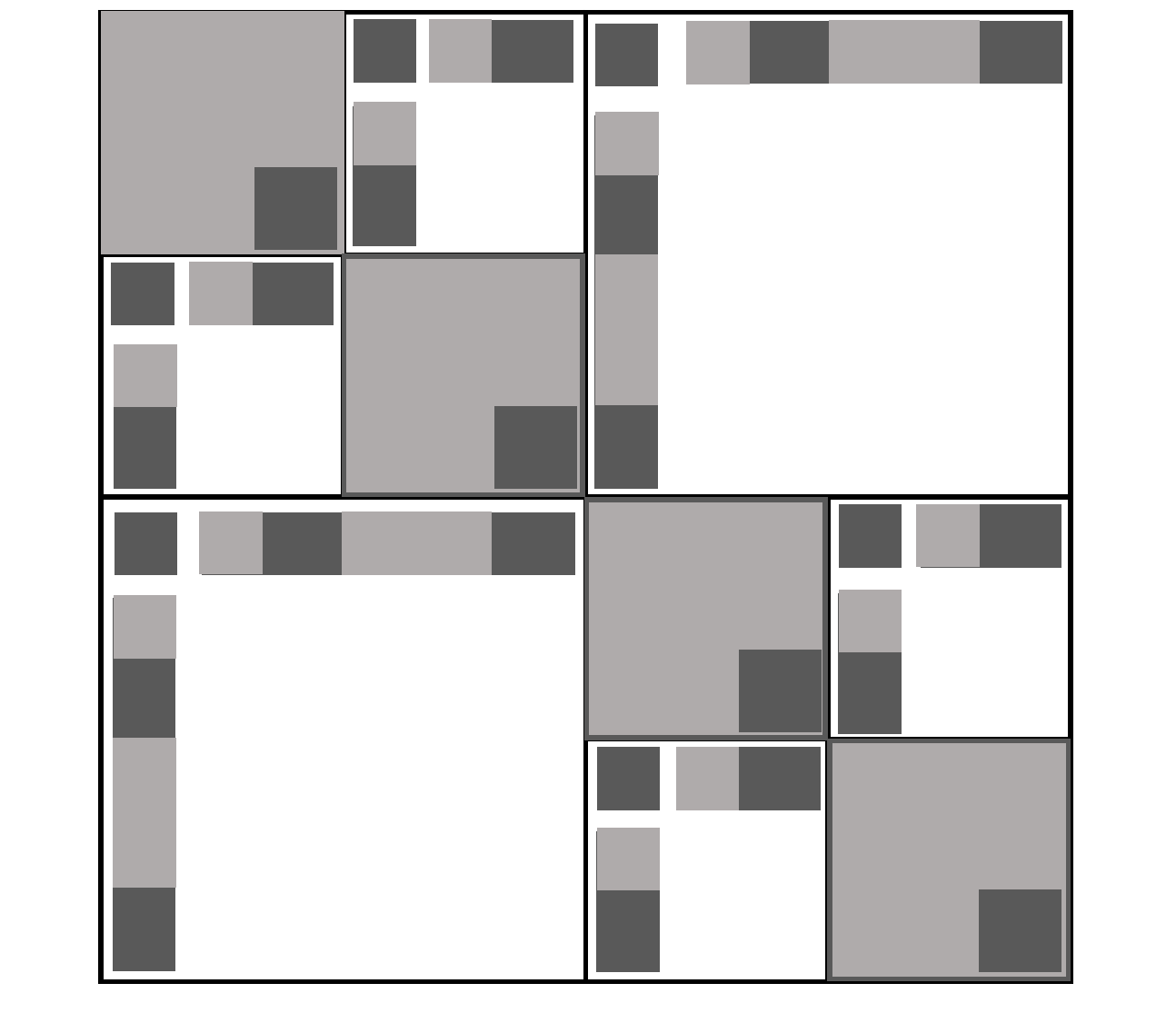}
 \caption{skeletonization}
\end{subfigure}
\begin{subfigure}{0.27\textwidth}
 \centering
 \includegraphics[width=\textwidth,trim=4.25cm 8.4cm 4.2cm 8.1cm,clip]{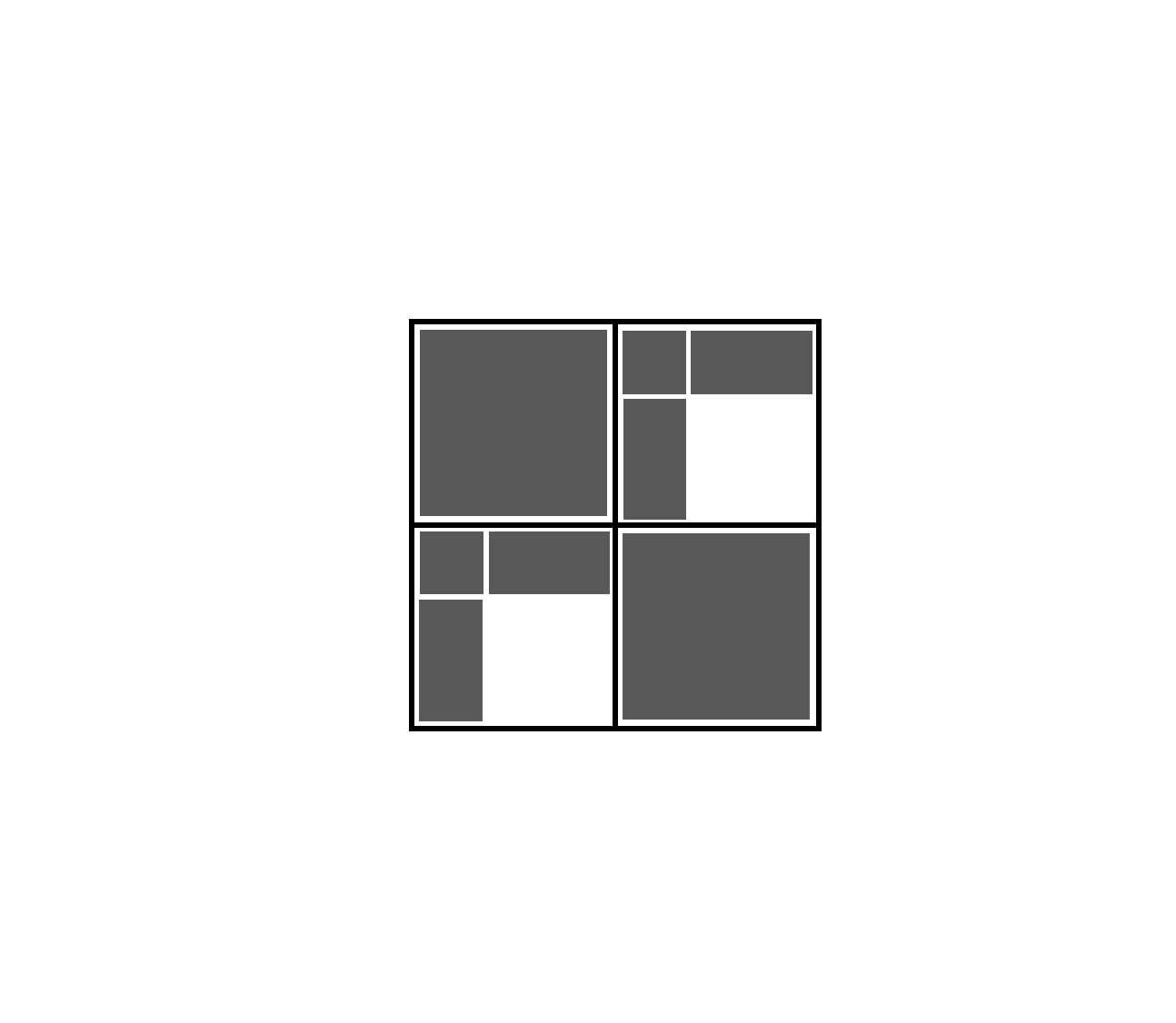}
 \caption{merge children}
\end{subfigure}\\
\begin{subfigure}{0.9\textwidth}
 \centering
 \includegraphics[width=\textwidth,trim=3.5cm 12.75cm 2.55cm 11.0cm,clip]{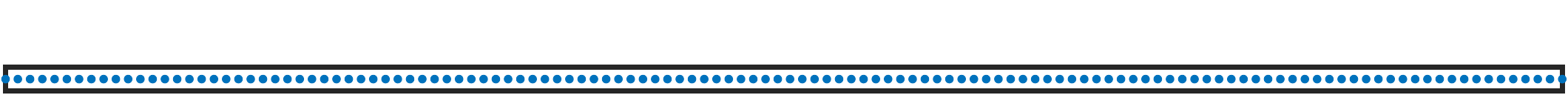}
 \caption{Start factorization}
\end{subfigure}\\
\begin{subfigure}{0.9\textwidth}
 \centering
 \includegraphics[width=\textwidth,trim=3.5cm 12.75cm 2.55cm 11.8cm,clip]{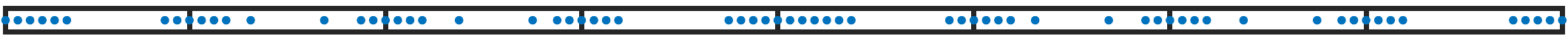}
 \caption{$\ell=3$ skeletonization}
\end{subfigure}\\
\begin{subfigure}{0.9\textwidth}
 \centering
 \includegraphics[width=\textwidth,trim=3.5cm 12.75cm 2.55cm 11.8cm,clip]{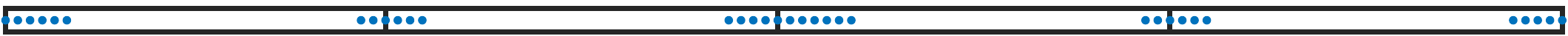}
 \caption{$\ell=2$ skeletonization}
\end{subfigure}\\
\begin{subfigure}{0.9\textwidth}
 \centering
 \includegraphics[width=\textwidth,trim=3.5cm 12.75cm 2.55cm 11.8cm,clip]{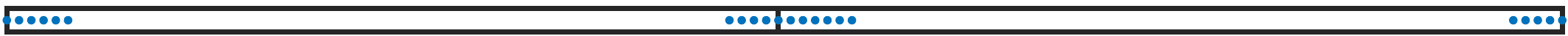}
 \caption{$\ell=1$  skeletonization}
\end{subfigure}
\caption{Top row, from left to right: original matrix $A$, skeletonization step zeroing the light
  gray shade DOFs of the matrix and maintaining the skeleton DOFs in dark gray, merging the
  skeletons with common parent node.  Rows 2 to 5 show the active DOFs on the 1D domain $\Omega$ after
  each step of the factorization, on three successive levels of skeletonization on the tree
  $\Tree_{\Omega}$.}
\label{fig:RSS}
\end{figure}

At the end of this process, we have the resulting matrix 
\begin{align} \label{eq1d_A}
 \Hm_L \approx R_{1}^* \cdots R_L^* \Hm R_L \cdots R_{1}
\end{align}
at the root level of the tree, being the identity everywhere except in the submatrix indexed by the
remaining active DOFs $\S_0$, which include DOFs near the boundary of the two subdomains $\Omega$ is
divided into at level $\ell=1$.  One can easily write an approximate factorization of $\Hm$ as
\begin{align} 
 \Hm \approx (R_{L}^*)^{-1} \cdots (R_{1}^*)^{-1}\Hm_{L}R_{1}^{-1} \cdots R_{L}^{-1},
 \label{eqA_1d}
\end{align}
where the matrices of the form $R_{\ell}$ can be easily inverted since they are
triangular up to a permutation.  By inverting the factors on can efficiently get an approximate
factorization of the inverse of $\Hm$ as
\begin{align}
  \Ga \equiv R_{L} \cdots R_{1}\Hm_{L}^{-1}R_{1}^{*} \cdots R_{L}^{*} \approx \Hm^{-1}=(\Sa)^{-1} \approx (K^*K)^{-1}
  \label{eqAinv_1d}
\end{align} 
that can be applied in linear time and whose error is controlled by the tolerance $\epsilon$ used in
the skeletonization steps.  In Figure \ref{fig:RSS} we show the active DOFs after each level of the tree following the inversion algorithm, and how the skeletonization and merging steps look like for the HODLR representation.

Note that in the skeletonization step we first compute an interpolative decomposition and then we apply Gaussian elimination on the leading block. Alternatively, as shown in \cite{XiaHSS10} we could first compute the Cholesky decomposition of the diagonal block $\A_{\I\I} = L_{\I}L_{\I}^*$ and scale the matrix
\begin{align}
    \begin{bmatrix}
           L_{\I}^{-1} & \\
           & I 
    \end{bmatrix} \begin{bmatrix}
           \A_{\I\I} & \A_{\I\R}\\
           \A_{\I\R}& \A_{\R\R} 
    \end{bmatrix} \begin{bmatrix}
           (L_{\I}^*)^{-1} & \\
           & I 
    \end{bmatrix} = \begin{bmatrix}
           I & L_{\I}^{-1}\A_{\I\R}\\
           \A_{\I\R}(L_{\I}^*)^{-1}& \A_{\R\R} 
    \end{bmatrix}
    \label{eq:ID_1d_resc}
    \end{align}
so that the leading block becomes the identity matrix. Next, one can perform interpolative decomposition on $\bar{U}_{\ell;\I}(L_{\I}^{*})^{-1}$. In this case, after performing ID there is no need to perform Gaussian elimination \eqref{eq_sk_3} because the leading block is already identity. 
This may be advantageous in the case of ill-conditioned matrices, leading to improved accuracy of the inverse approximation \cite{Xia10}.

\subsection{2D case}
In the multidimensional case, we have constructed an $\mathcal{H}$-matrix approximation $\Sa\approx
S\equiv K^*K$. Denoting $\Sa$ by $\A$ in this subsection, we can invert such representation in a
similar way to that of the hierarchical interpolative factorization \cite{HIFIE}, but exploiting the
$\mathcal{H}$-matrix representation.

Similarly to the one-dimensional case, the factorization proceeds from the leaves level $\ell=L$ of
the tree $\Tree_{\Omega}$ to level $\ell=1$ and performs the following three steps at each level
$\ell$:

\begin{itemize}
\item \textbf{Skeletonization of cells at level $\ell$.}  There are $q_{\ell} = 2^\ell\times 2^\ell$ square
  cells at level $\ell$. For each cell $I$ we let $\I$ be the set of indices for the active DOFs in
  the cell. Then, for each set of indices $\I$, we define the matrix
  \begin{align}
    \bar{U}_{\ell;\I} =     \begin{bmatrix}
          \Hm_{\Bi\I}\\ \Ub_{\ell;\I}^*\\
           \vdots\\
           \Ub_{1;\I}^*
    \end{bmatrix},
    \label{eq:2d_id_face}
  \end{align}
  with $\Hm_{\ell;\Bi\I}$ the block of matrix $\Hm_{\ell}$ indexed by DOFs in $\I$ and $\Bi$, where
  $\Bi$ includes the active DOFs of the neighboring cells $J\in NL(I)$ in level $\ell$ of the tree,
  i.e. the near field of $I$. Denote $I_{\ell}$ the node containing the DOFs $\I$ at level $\ell$ of $\T_{\Omega}$ and $m$ the number of nodes in $IL(I_{\ell})$. The submatrices of the form $\Ub_{\ell;I}$ contain the low rank bases at level $\ell$ for the low-rank approximation of subblocks $\A_{\J_{\ell,i}\I_{\ell}} = \U_{\ell;\J_{\ell,i}\I_{\ell}}\B_{\ell;\J_{\ell,i}\I_{\ell}}\U_{\ell;\J_{\ell,i}\I_{\ell}}^*$ for all nodes $J_{\ell,i}\in IL(I_{\ell})$,
  \begin{align}
      \Ub_{\ell;I} = \begin{bmatrix}
           B_{\ell;\J_{\ell,1}\I_{\ell}}\U_{\ell;\J_{\ell,1}\I}^*\\ B_{\ell;\J_{\ell,2}\I_{\ell}}\U_{\ell;\J_{\ell,2}\I}^*\\
           \vdots\\  B_{\ell;\J_{\ell,m}\I_{\ell}}\U_{\ell;\J_{\ell,m}\I}^*
    \end{bmatrix},
    \label{eq:2d_bases}
  \end{align}
  with $\U_{\ell;\J_{\ell,i}\I}$ containing only the columns indexed by $\I$ from $\U_{\ell;\J_{\ell,i}\I_{\ell}}$.
 
  Similar to the 1D case, we perform interpolative
  decomposition on $\bar{U}_{\ell,\I}$ \eqref{eq_id} to obtain the corresponding interpolation matrix $T_{\I}$. Now one can apply the zeroing matrix $\Z_{\I}$ \eqref{eq_sk} followed by Gaussian elimination to eliminate
  the redundant DOFs $\Ir$ as in \eqref{eq_sk_3},
  \begin{align} \label{eq_sk_2d}
\E_{\Ir}^{*}\Z_{\I}^*\A_{\ell} \Z_{\I} \E_{\Ir} \approx \begin{bmatrix}
  I&&\\
  &\He_{\Is\Is}&\Hm_{\R\Is}^*\\
  &\Hm_{\R\Is}&\Hm_{\R\R}\\
\end{bmatrix}.
\end{align}
  
  Notice we are computing $\T_{\I}$ using $\bar{U}_{\ell;\I}$ instead of the subblock $\Hm_{\R\I}$. This is equivalent
  but essential to reduce the numerical complexity of the algorithm, since $\bar{U}_{\ell;\I}$ has
  many less rows than $\Hm_{\R\I}$. In particular one can write $\Hm_{\R\I}$ as
  \begin{align}
    \Hm_{\R\I} =
    \begin{bmatrix}
      I&&&\\
      &V_{\ell;\I}&&\\
      && \ddots&\\
      &&&V_{1;\I}
    \end{bmatrix}
    \begin{bmatrix}
      \Hm_{\Bi\I} \\
      \Ub_{\ell;\I}^*\\
      \vdots\\
      \Ub_{1;\I}^
      *
    \end{bmatrix} \equiv  Q_{\ell;\I}\bar{U}_{\ell;\I}
  \end{align}
  with $Q_{\ell;\I}$ having orthogonal columns and $V_{\ell;\I}$ containing the low-rank bases of all $m$ nodes in $IL(I_{\ell})$,
  \begin{align}
      V_{\ell;\I} =
    \begin{bmatrix}
      U_{\ell;\J_1\I_{\ell}}&&&\\
      &U_{\ell;\J_2\I_{\ell}}&&\\
      && \ddots&\\
      &&&U_{\ell;\J_m\I_{\ell}}
    \end{bmatrix}
    \label{eq:column_bases}
  \end{align}
  
  with $U_{\ell;\J_i\I_{\ell}}$ the block $U_{\ell;\J_{\ell,i}\I_{\ell}}$ (from the low-rank approximation $A_{\J_{\ell,i}\I_{\ell}} = U_{\ell;\J_{\ell,i}\I_{\ell}}B_{\ell;\J_{\ell,i}\I_{\ell}}U_{\ell;\J_{\ell,i}\I_{\ell}}^*$ in the $\mathcal{H}$-matrix representation) restricted to the rows indexed by the active DOFs of the node $J_{\ell,i}\in IL(I)$. Performing ID on $\bar{U}_{\ell;\I}$ gives an error
  \begin{align}
    \begin{split}
      \|E_{\I}\| =& \|\Hm_{\R\Ir}-\Hm_{\R\Is}\T_{\I}\| = \|Q_{\ell,\I} \bar{U}_{\ell;\Ir}-Q_{\ell,\I}\bar{U}_{\ell;\Is}\T_{\I}\| =\\ &\|\bar{U}_{\ell;\Ir}-\bar{U}_{\ell;\Is}\T_{\I}\| = O(\epsilon\|\bar{U}_{\ell,\I}\|)= O(\epsilon\|\Hm_{\R\I} \|),
    \end{split}
  \end{align}
  of the same order of magnitude that if ID had been performed on $\Hm_{\R\I}$ directly, since
  $Q_{\ell,\I}$ has orthogonal columns, and therefore $\|\Hm_{\R\I} \| =
  \|\bar{U}_{\ell,\I}\|$.
    
  Let $\{ \I_{\ell, i} \}_{i = 1}^{q_\ell}$ be the collection of disjoint index sets corresponding to active DOFs of each of the $q_{\ell}$ cells
  at level $\ell$. Performing skeletonization over all cells at level $\ell$ results in
  
    \begin{align} 
    \He_{\ell}\approx R_{\ell}^{*}\Hm_{\ell}R_{\ell}, \qquad
    R_{\ell}=\prod_{i=1}^{q_\ell}\Z_{\I_{\ell,i}}E_{\Ir_{\ell,i}},
  \end{align}
  
  where we have decoupled all
  the redundant DOFs in the square cells. Now, the remaining active DOFs in $\He_{\ell}$ are
  $\cup_{i=1}^{q_\ell}\Is_i$, and the DOFs $\cup_{i=1}^{q_\ell}\Ir_i$ have been decoupled,
  i.e. their corresponding matrix block in $\He_{\ell}$ is the identity matrix.
    
  Analogously to the 1D case, for each cell at a level $\ell$, with active DOFs indexed by $\I$, one
  can scale the off-diagonal low-rank bases with the Cholesky factor of $\A_{\ell;\I\I}$ before
  computing the interpolative decomposition for improved accuracy of the factorization for ill-conditioned matrices.
    
\item \textbf{Skeletonization around the edges at level $\ell$.}  After performing cell
  skeletonization, most skeletons concentrate around the edges. This happens because of the
  high-rank interaction between the DOFs of adjacent cells located near a common edge. However, if
  these nodes are considered all together, then their interaction with the rest of the matrix can be
  assumed to be low-rank and we can reduce the number of skeleton DOFs around each edge.
    
  In order to take advantage of this behavior, we start by assigning every active DOF to the
  nearest edge at level $\ell$. Each edge will have DOFs assigned to it from two adjacent cells.
  For each edge, let $\I$ be the active DOFs assigned to it from adjacent cells. Then, we define the
  matrix
  \begin{align}
    \label{ID_2d_bases}
    \bar{U}_{\ell;\I} =     \begin{bmatrix}
      \He_{\Bi\I_1} &\He_{\Bi\I_2}\\ \Ub_{\ell;\I_1}^*& \\
      \vdots&\\
      \Ub_{1;\I_1}^*&\\
      & \Ub_{\ell;\I_2}^*\\
      &\vdots\\
      &\Ub_{1;\I_2}^*\\
    \end{bmatrix}
  \end{align}
  with $\I_1$ and $\I_2$ the indices of $\I$ corresponding to each of the two cells adjacent to the
  edge and subblocks $\Ub_{\ell,\I_i}$ defined as in \eqref{eq:2d_bases}. Here, $\Bi$ contains the
  DOFs in the union of the near field and the interior DOFs of the cells adjacent to the edge we are
  sparsifying.  Now we can compute the interpolative decomposition of $\bar{U}_{\ell;\I}$ to
  construct the zeroing matrix $\Z_{\I}$ \eqref{eq_sk} and decouple the redundant DOFs $\Ir$ as in
  \eqref{eq_sk_3} with Gaussian elimination resulting in
  \begin{align} \label{eq_edge_sk_2d}
    \E_{\Ir}^{*}\Z_{\I}^*\He_{\ell} \Z_{\I} \E_{\Ir} \approx \begin{bmatrix}
      I&&\\
      &\bar{\He}_{\Is\Is}&\He_{\R\Is}^*\\
      &\He_{\R\Is}&\He_{\R\R}\\
    \end{bmatrix}.
  \end{align}
  
  One can decompose $\He_{\R\I}$ as
  \begin{align}
    \He_{\R\I} = 
    \begin{bmatrix}
      I&&&&&\\
      &V_{\ell;\I_1}&&&V_{\ell;\I_2}&&\\
      && \ddots&&& \ddots&\\
      &&&V_{1;\I_1}&&&V_{1;\I_2} \\
    \end{bmatrix}
    \bar{U}_{\ell,\I}  \equiv  Q_{\ell,\I}\bar{U}_{\ell,\I}
  \end{align}
  where $V_{\ell;\I_1}$ and $V_{\ell;\I_2}$ are defined as in \eqref{eq:column_bases}. Analogously to cell skeletonization, the error committed with ID performed on $\bar{U}_{\ell;\I}$ is 
  \begin{align}
    \begin{split}
      \|E_{\I}\| =& \|\Hm_{\R\Ir}-\Hm_{\R\Is}\T_{\I}\| = \|Q_{\ell,\I} \bar{U}_{\ell;\Ir}-Q_{\ell,\I}\bar{U}_{\ell;\Is}\T_{\I}\| \leq\\ &\sqrt{2}\|\bar{U}_{\ell;\Ir}-\bar{U}_{\ell;\Is}\T_{\I}\| = O(\epsilon\|\bar{U}_{\ell,\I}\|),
    \end{split}
  \end{align}
  where we have used the fact $\|Q_{\ell,\I}\|\leq2$. Since $Q_{\ell,\I}$ is a tall skinny matrix
  composed of two tall-skinny matrices ($Q_{\ell,\I} = [Q_1, Q_2]$) with orthogonal columns,
  concatenated horizontally, it is trivial to bound its norm by
  \begin{align}
    \begin{split}
      \|(\begin{matrix}Q_1 & Q_2\end{matrix})\| =& \max_{\|(\begin{matrix}v_1 & v_2\end{matrix})\|=1}\|(\begin{matrix}Q_1 & Q_2\end{matrix})(\begin{matrix}v_1 & v_2\end{matrix})^T\|\leq\max_{\|(\begin{matrix}v_1 & v_2\end{matrix})\|=1}\|Q_1v_1^T+Q_2v_2^T\|\leq \\
            &\max_{\|(\begin{matrix}v_1 & v_2\end{matrix})\|=1}\|Q_1v_1^T\|+\|Q_2v_2^T\|= \max_{\|(\begin{matrix}v_1 & v_2\end{matrix})\|=1}\|v_1\|+\|v_2\|=\sqrt{2}.
    \end{split}
  \end{align}
  
  Let $s_{\ell}$ be the number of edges and $\{ \I_{\ell, i} \}_{i = 1}^{s_\ell}$ the collection of corresponding disjoint index sets
  at level $\ell$, corresponding to each edge $i$. Performing skeletonization on all the index sets
  gives
  \begin{align} 
    \Hm_{\ell-1}\approx\bar{R}_{\ell}^{T}\He_{\ell}\bar{R}_{\ell}, \qquad
    \bar{R}_{\ell}=\prod_{i=1}^{s_\ell}\Z_{\I_{\ell,i}}\E_{\Ir_{\ell,i}},
  \end{align}
  where the remaining active DOFs in $\Hm_{\ell-1}$ are $\S_{\ell-1} =
  \cup_{i=1}^{s_\ell}\Is_{\ell,i}$. The index sets $\cup_{i=1}^{s_\ell}\Ir_{\ell,i}$ have been
  decoupled, i.e. their corresponding matrix block in $\Hm_{\ell-1}$ is the identity matrix.
    
\item \textbf{Merge child blocks.}  After the skeletonization step on all edges, we merge the active DOFs of
  sibling nodes and move to the next level in the tree, i.e. to the parent node at level $\ell-1$
  obtained from merging 4 adjacent square cells into a bigger square cell.
\end{itemize}

\begin{figure}
\centering \captionsetup[subfigure]{labelformat=empty}
\begin{subfigure}{0.25\textwidth}
 \centering
 \includegraphics[width=\textwidth,trim=3.20cm 1.65cm 2.75cm 1.15cm,clip]{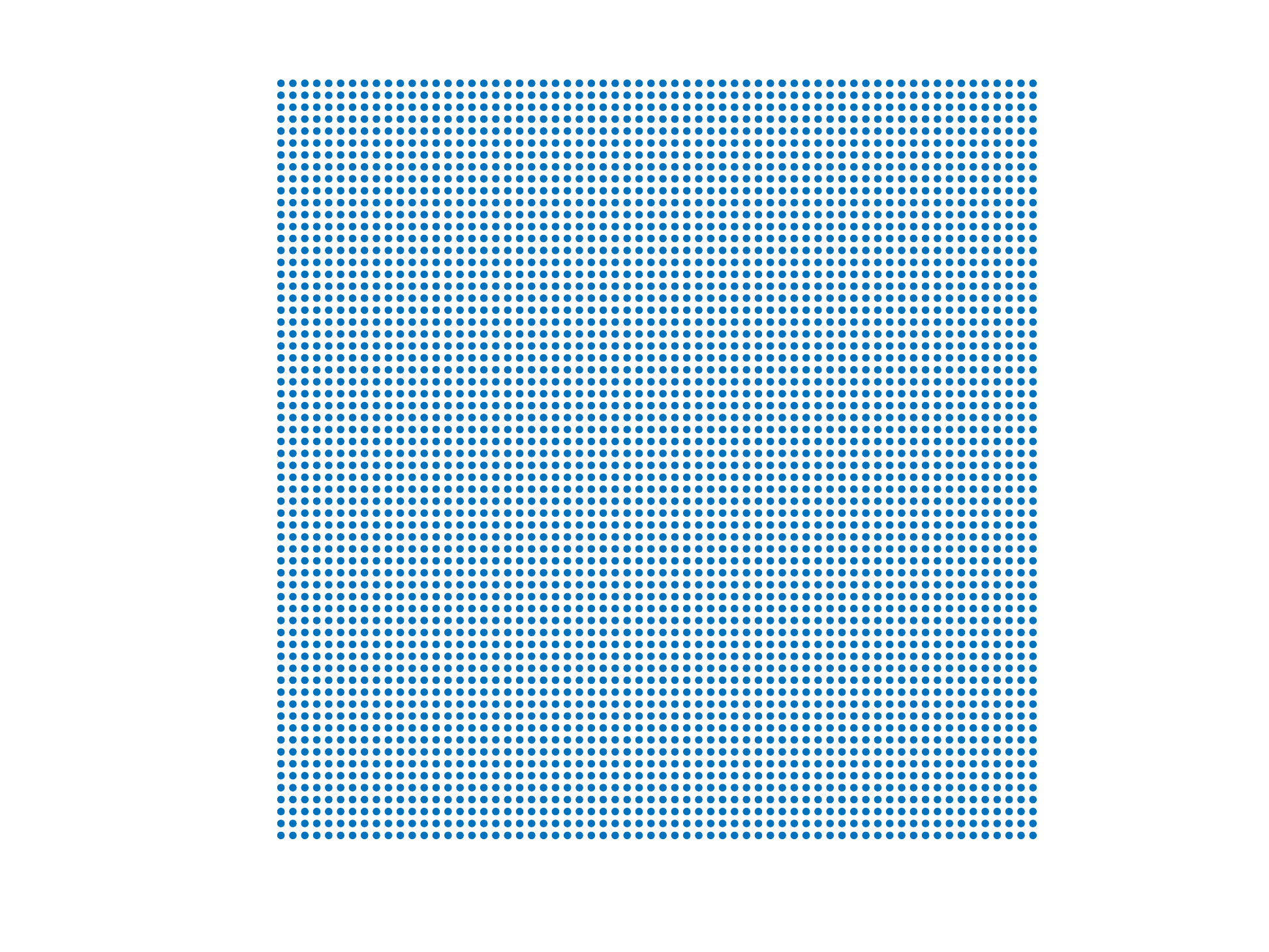}
 \caption{factorization start}
\end{subfigure}
\begin{subfigure}{0.25\textwidth}
 \centering
 \includegraphics[width=\textwidth,trim=3.20cm 1.65cm 2.75cm 1.15cm,clip]{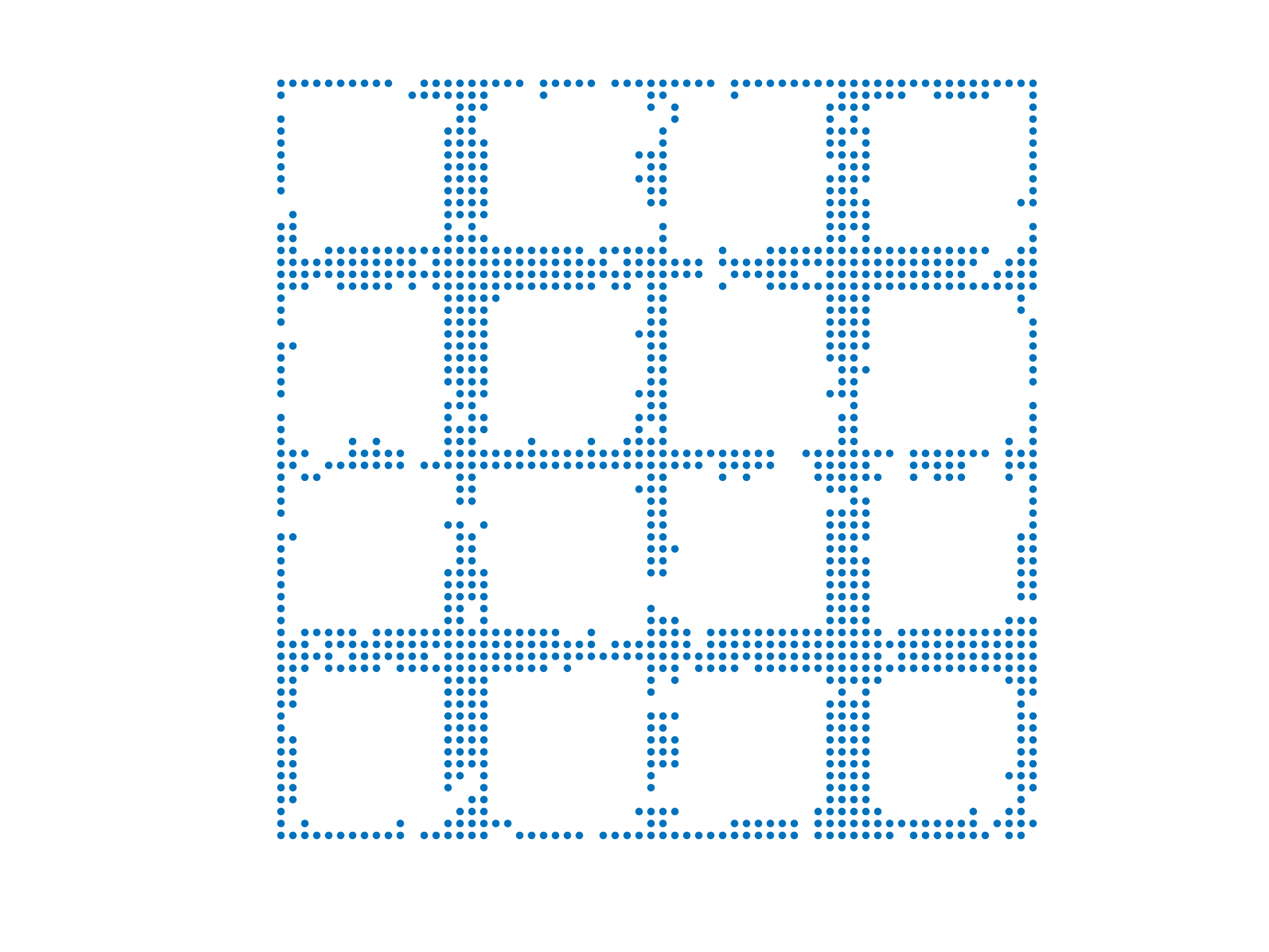}
 \caption{face skeletonization}
\end{subfigure}
\begin{subfigure}{0.25\textwidth}
 \centering
 \includegraphics[width=\textwidth,trim=3.20cm 1.65cm 2.75cm 1.15cm,clip]{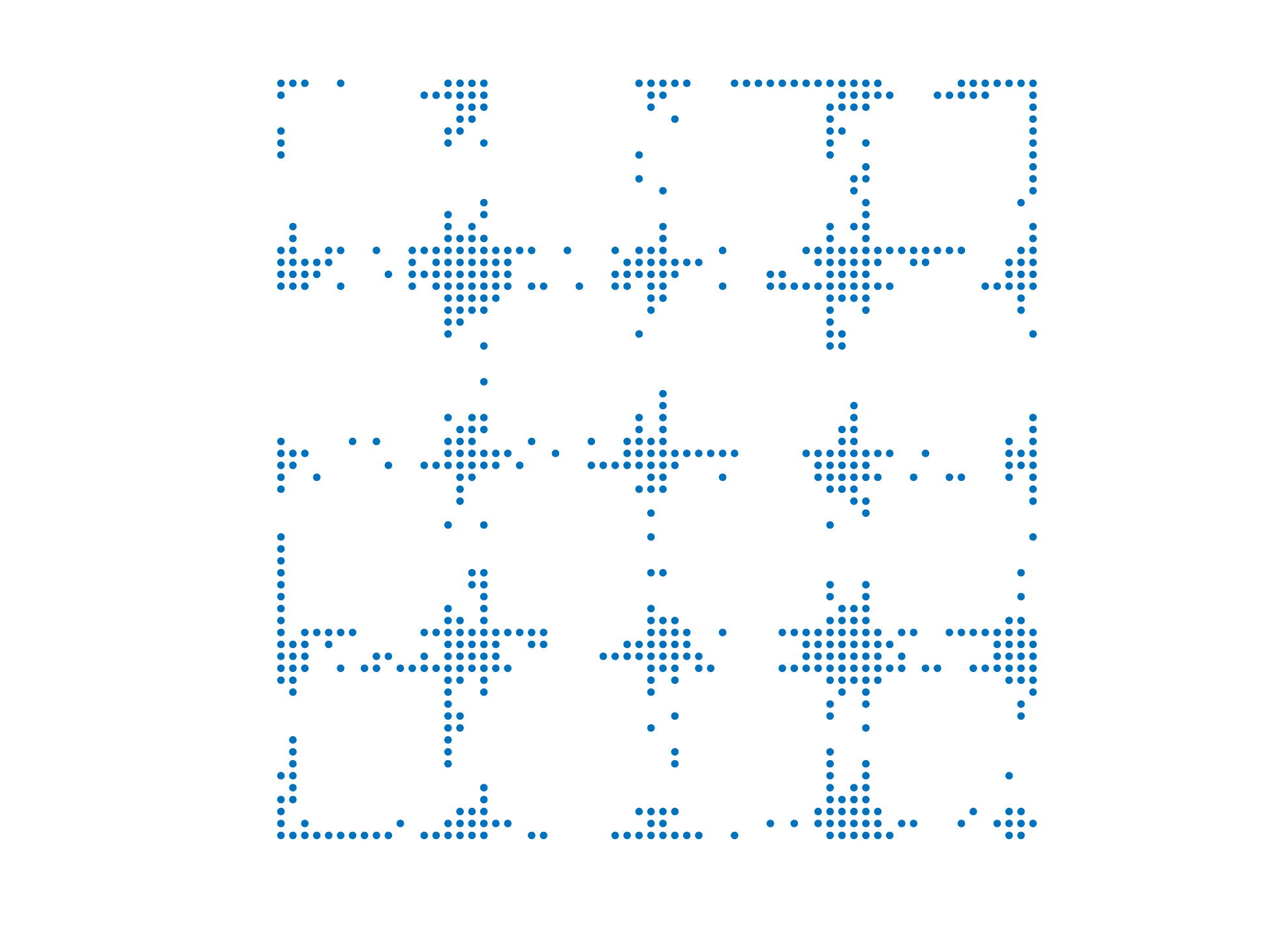}
 \caption{edge skeletonization}
\end{subfigure}\\
\begin{subfigure}{0.25\textwidth}
 \centering
 \includegraphics[width=\textwidth,trim=3.20cm 1.65cm 2.75cm 1.15cm,clip]{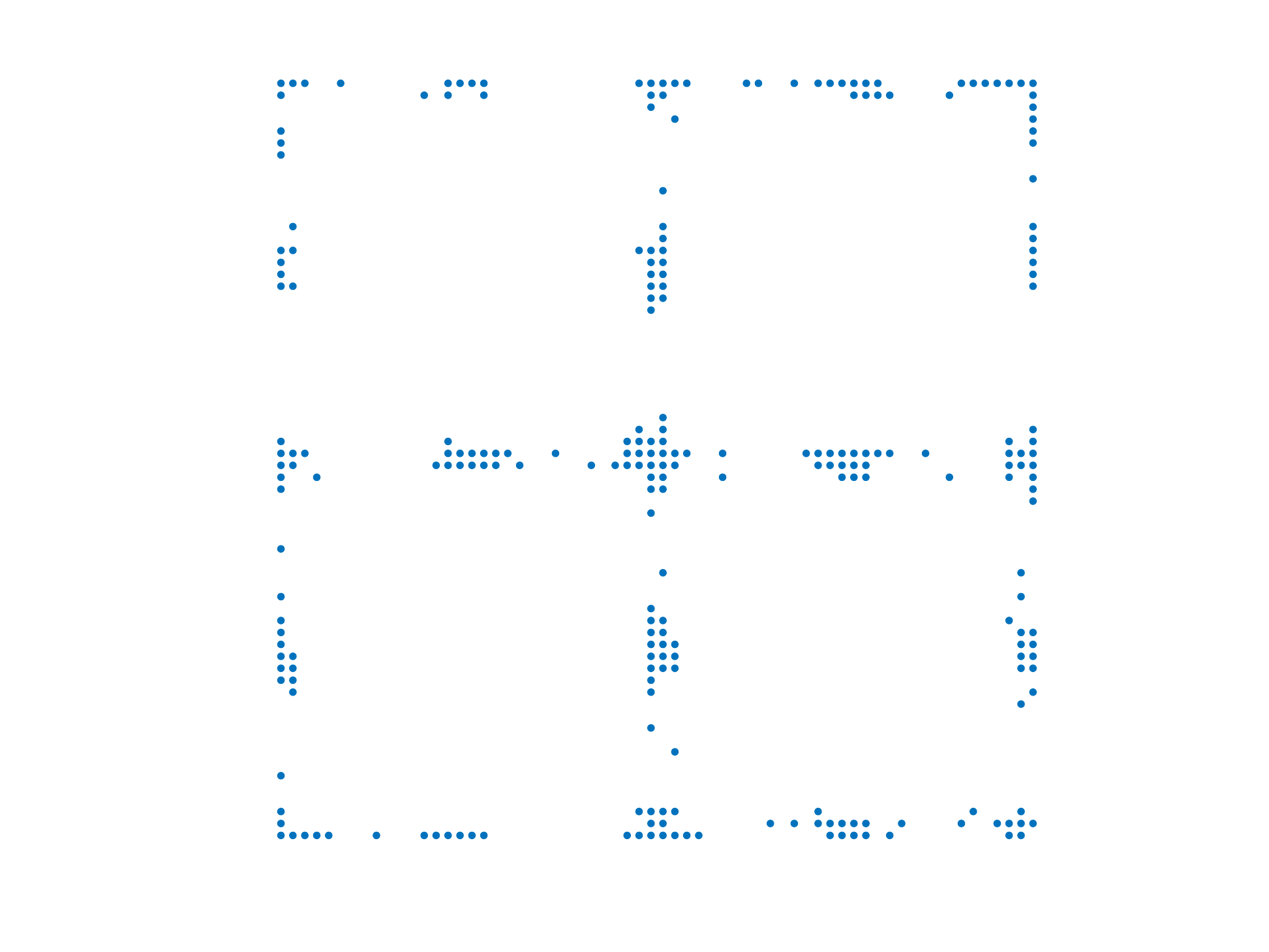}
 \caption{face skeletonization}
\end{subfigure}
\begin{subfigure}{0.25\textwidth}
 \centering
 \includegraphics[width=\textwidth,trim=3.20cm 1.65cm 2.75cm 1.15cm,clip]{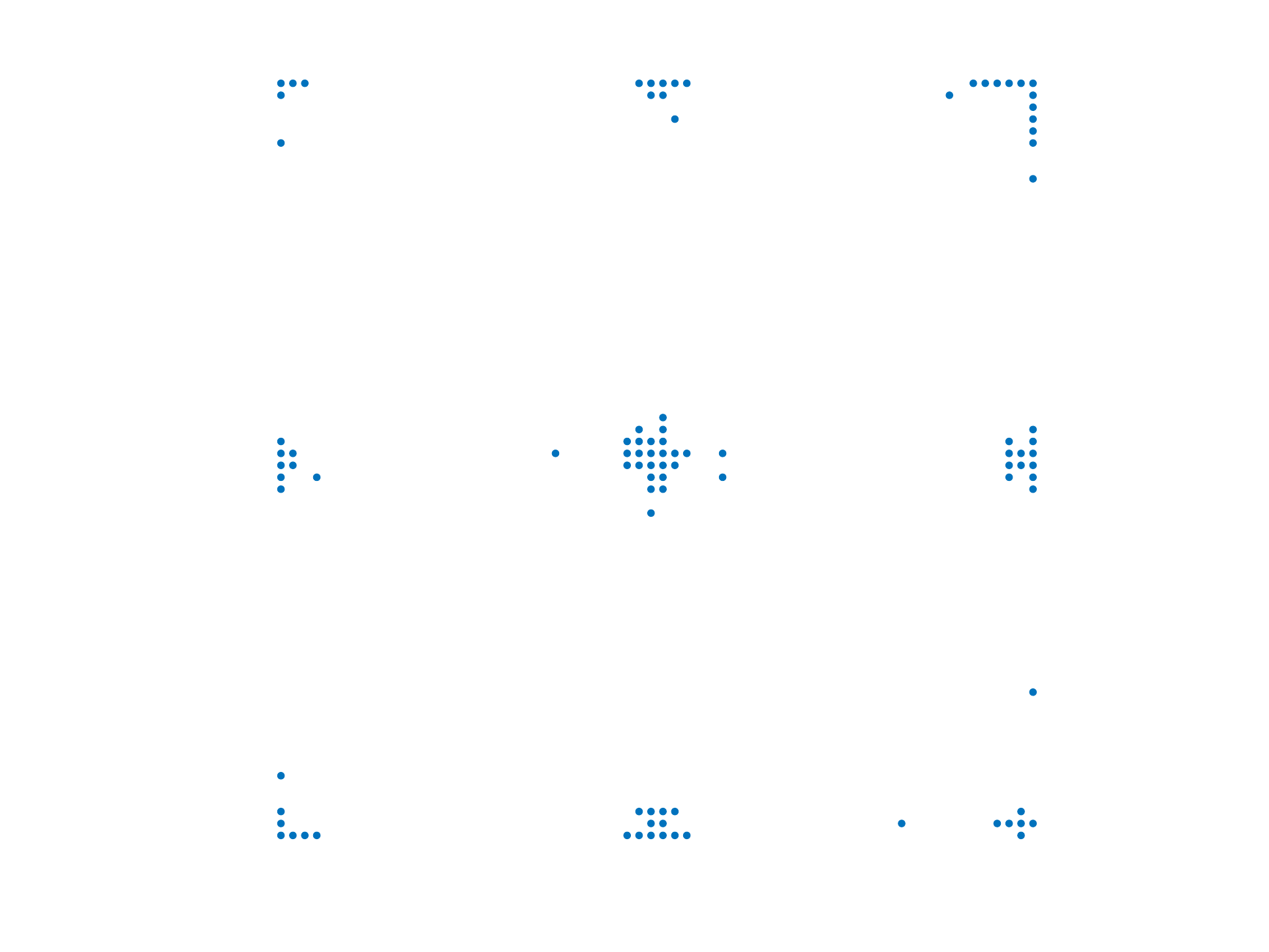}
 \caption{edge skeletonization}
\end{subfigure}
\caption{Active DOFs on the 2D square domain $\Omega$ for each step of the factorization, on two successive levels of the tree $\Tree_{\Omega}$}
\label{fig:HIFDE}
\end{figure}

Similarly to the 1D case, at the end of the factorization we can approximate the inverse of the $\mathcal{H}$-matrix as 
\begin{align}
  \Ga = R_{L} \bar{R}_{L}R_{L-1} \bar{R}_{L-1} \cdots
  R_{1}\bar{R}_{1}\Hm_{L}^{-1}\bar{R}_{1}^{*}R_{1}^{*} \cdots
  \bar{R}_{L-1}^{*}R_{L-1}^{*}\bar{R}_{L}^{*}R_{L}^{*}
  \approx \Hm^{-1}=(\Sa)^{-1} \approx (K^*K)^{-1},
  \label{eqAinv_2d}
\end{align}
with $R_{\ell}$ and $\bar{R}_{\ell}$ the cell and edge skeletonization matrices at level $\ell$
respectively. The matrix $A_L$ defined at the root level can be easily inverted since it has few
remaining active DOFs.  Figure \ref{fig:HIFDE} illustrates the active DOFs concentrating around the
edges and the vertices after each cell and edge skeletonization steps respectively, for two levels
of the tree.

\section{Approximate inverse factorization for FIOs}\label{sec:Alg}

We now proceed to describe the algorithm for inverting the discrete FIO.  Given that the discrete
FIO matrix $\K$ is square and invertible, we can obtain its inverse by computing the pseudo-inverse
\begin{align}
  \K^{-1} = (\Kt\K)^{-1}\Kt.
  \label{eq:pseudoinv}
\end{align}
First, we approximate $\Kt$ and $\K$ using the butterfly factorization $\Ka \approx \K$. Since the
Hermitian matrix $\Kt\K$ has low-rank blocks for well separated clusters, we can use HODLR and
$\mathcal{H}$-matrix representations in 1D and 2D respectively.

The HODLR (or $\mathcal{H}$-matrix) approximation $\Sa$ of $S\equiv \Kt\K$ is then
computed using the peeling algorithm. Fast matrix-vector multiplications with $\K^*\K$ needed in
the peeling algorithm are obtained by multiplying with $\Ka$ and $\Ka^*$, which can be done in quasi-linear time.

Next, we invert the hierarchical matrix approximation $\Sa$ of $\Kt\K$ using the inversion algorithm described in section \ref{sec:HIF}, obtaining an approximate inverse $\Ga \approx (\Kt\K)^{-1}$.  Finally, we can
construct the approximate factorization of the inverse FIO by approximating the terms $(\Kt\K)^{-1}$
and $\Kt$ in \eqref{eq:pseudoinv} with $\Ga$ and $\Ka^*$ respectively, i.e. $ \K^{-1}\approx \Ga\Ka^*$. The algorithm is
summarized in Algorithm \ref{alg}.

One can use the above factorization to get an approximate inverse, but for higher accuracy, $\Ga\Ka^*$
can be used as a preconditioner.  If we want to solve a non-symmetric system of equations of the
form $\K F f = u$ \eqref{eq6}, we can first solve the system
\begin{align}\K \hat{f} = u
\label{eq:system}\end{align}
and then apply inverse FFT to $f$. In order to solve \eqref{eq:system} 
efficiently, we can first approximate the FIO matrix by its butterfly factorization $\Ka\approx \K$
so that matrix-vector multiplication can be performed in quasi-linear time. Next, we precondition
the system with $\Ga\Ka^*$ and solve the approximate system
\begin{align} \Ga\Ka^*\Ka \hat{f}= \Ga\Ka^*u \end{align}
for which we can run conjugate gradient \cite{Hestenes} on the hermitian matrix $\Ka^*\Ka$
preconditioned with $\Ga$, $\text{pcg}(\Ka^*\Ka,\text{ }\Ka^*u,\text{ } \Ga)$.

\begin{algorithm}[htb]
  Given: FIO phase function $\Phi(x,\xi)$, butterfly factorization tolerance $\epsilon_{BFF}$, peeling tolerance $\epsilon_{peel}$, HIF inversion tolerance $\epsilon_{HIF}$\\
  $//$ \texttt{Factor the FIO matrix $\K$ with butterfly factorization with tolerance} $\epsilon_{BFF}$\\
  $\Ka \leftarrow \text{BFF}(\Phi(x,\xi),\epsilon_{BFF})$\\
  $//$ \texttt{Construct Hierarchical matrix representation of $\Kt\K$ with matrix peeling using fast matrix-vector multiplications with} $F$ \texttt{and tolerance} $\epsilon_{peel}$\\
  $\Sa \leftarrow \text{peeling}(\Ka,\epsilon_{peel})$\\
  $//$ \texttt{Invert} $\Sa$ \texttt{with the hierarchical interpolative factorization}\\
  $\Ga \leftarrow \text{HIF}(\Sa,\epsilon_{HIF})$\\
  $//$ \texttt{Construct the inverse FIO}\\
  $\K^{-1}\approx \Ga \Ka^*$

  \caption{Inverse factorization of multidimensional FIOs}
  \label{alg}
\end{algorithm}

\subsection{Computational complexity}\label{complexity}
The butterfly factorization $\Ka$ can be computed in $O(N^{3/2})$ time complexity, and its direct
and transpose application to a vector $v$, $\Ka v$ and $\Ka^* v$, can be computed in quasi-linear time
$O(N\log N)$.

The HODLR and $\mathcal{H}$-matrix approximation $\Sa$ of $\Kt\K$ can be constructed in
$O(N\log^2N)$, since $\Ka^*\Ka$ can be applied in $O(N\log N)$ time to any arbitrary vector $v$ and
we need to apply it to $O(\log N)$ proving matrices with size $N\times(k+c)$.

Finally, $\Sa$ can be inverted in $O(N\log^2N)$ time  since during the factorization we need to perform interpolative decomposition on matrices resulting from concatenation of $O(\log N)$ submatrices. Such factorization can be applied in quasi-linear time $O(N\log\log N)$
according to \cite{HIFIE}.

Therefore, the computational complexity for computing a factorization for the inverse FIO matrix
$\K^{-1}$ is $O(N^{3/2})$ and its application to a vector is $O(N\log N)$. Table
\ref{table:complexity} summarizes the asymptotic complexity of each step in the construction and
application of the inverse factorization of discrete FIOs.

\begin{table}[htb]
\centering

\begin{tabular}{c|c}
\toprule
Operation & Complexity\\
\midrule
Construct $\Ka$ & $O(N^{3/2})$\\
Apply $\Ka$ and $\Ka^*$ to a vector & $O(N \log N)$\\
Construct $\mathcal{H}$-matrix $\Sa$ & $O(N \log^2 N)$\\
Construct $\Ga$ & $O(N \log^2N)$\\
Apply $\Ga$ to a vector & $O(N \log\log N)$ \\
Apply $\Ga\Ka^*\approx \K^{-1}$ to a vector & $O(N \log N)$\\
\bottomrule

\end{tabular}
\caption{Complexity estimates for different steps in the FIO inverse factorization}
\label{table:complexity}
\end{table}

\section{Numerical results}\label{sec:NR}
In this section, we demonstrate the effectiveness of the method for 1D and 2D FIOs. Here, we have
considered a typical setting with periodic functions $a(x,\xi)$ and $\phi(x,\xi)$ on $x$, with
$\Phi(x,\xi) = x\cdot\xi + \phi(x,\xi)$. This allows one to embed the problem in a periodic cell,
chosen here to be $[0,1)^d$.

The factorization obtained to approximate the inverse operator is also used as a preconditioner to
solve a system of equations of the form $\Ka x = b$. The performance of the preconditioner is
compared with $\Ka^*$, which can also be used as a preconditioner. Notice that
without a preconditioner any iterative solver will take too many iterations to converge due to the
fact that eigenvalues are not clustered.  In particular, we compare the number of conjugate gradient (CG) iterations of the system
\begin{align}
  \Ka^*\Ka x= \Ka^*b, \label{eq_NR_1}
\end{align}
preconditioned with and without $\Ga$.  The following quantities are used in the rest of the section
to evaluate the performance of the preconditioner:

\begin{itemize}
\item $N$: total number of DOFs in the problem;
\item $e_a$: the relative error of the butterfly approximation $\Ka$ of $\K$, defined as
  $\|\K-\Ka\|/\|\K\|$;
\item $\epsilon$: the tolerance used in the hierarchical matrix construction and inversion;
\item $e_s$: the relative error of the approximation $\Ga\Ka^*$ of $K^{-1}$, defined as
  $\|I-\Ga\Ka^{*} K\|\geq\|K^{-1}-\Ga\Ka^{*}\|/\|K^{-1}\|$;
\item $t_s$: total running time of CG to solve \eqref{eq_NR_1};
\item $n_i$: the number of CG iterations with  relative residual
  equal to $10^{-8}$.
\end{itemize}
The errors $e_a$ and $e_s$ are approximated using randomized power iteration \cite{Dixon,Kuczynski} to $10^{-2}$ relative precision. The butterfly factorization is computed using the code from \cite{BFcode}.

\paragraph{Example 1. 1D FIO with uniform amplitude.} We begin with an example of a one-dimensional discrete FIO of the form
\begin{align}
  \label{eq_fio_1d}
  u(x) = \int_{\Rb} a(x) e^{2\pi i \Phi(x,\xi)}\hat{f}(\xi)\text{d}\xi
\end{align}
with uniform amplitude $a(x,\xi) =1$ and phase function $\Phi(x,\xi)$ given by
\begin{align}
\Phi(x,\xi) = x\cdot\xi + c(x)|\xi|, \quad c(x) = (2 + \sin(2\pi x))/8.
\end{align}
Discretizing $x$ and $\xi$ on $[0,1)$ and $[-N/2,N/2)$ with $N$ points,
\begin{align}
    x_i =(i-1)/N \quad \text{and} \quad \xi_j =j-1-N/2,
\end{align}
leads to the discrete system $u = Kf$,
with $u_i = u(x_i)$, $f_j=\hat{f}(\xi_j)$ and $K_{ij}=e^{2\pi i \Phi(x_i,\xi_j)}$.

We compute the approximate inverse factorization with tolerance $\epsilon=10^{-3}$ and
$\epsilon=10^{-6}$, leading to errors $e_s$ of the order between $10^{-3}$ and $10^{-6}$ as reported
in Table \ref{1d_uniform}. The factorization $\Ga\Ka^*$ is also a much more efficient preconditioner
than the adjoint FIO matrix $\Ka^*$. For example, when $\epsilon=10^{-6}$ is used, CG converges in
just 2 iterations instead of 28, leading to a reduction in time by a factor of approximately 7.

\begin{table}[htb]
\centering
\begin{tabular}{c|c|cccc|cc}
\toprule
& $\Ka\approx \K $&\multicolumn{4}{|c|}{$\Ga\Ka^*\approx \K^{-1}$} & \multicolumn{2}{|c}{$\Ka^*\approx \Kt$}\\
\midrule
 $N$ & $e_a$ &$\epsilon$& $e_s$ & $n_i$ & $t_s$ & $n_i$ & $t_s$\\
\midrule
 \multirow{2}{*}{1024} & \multirow{2}{*}{$3.59e-7$} & {$1e-6$}& $4.32e-6$ &2 & $1.79e-2$ & \multirow{2}{*}{28} & \multirow{2}{*}{$1.17e-1$}\\
   &  & {$1e-3$}& $2.06e-3$ & 3 & $2.68e-2$ & &\\
 \midrule
 \multirow{2}{*}{4096} & \multirow{2}{*}{$4.77e-7$} &{$1e-6$}& $8.71e-6$ &2 & $1.15e-1$ & \multirow{2}{*}{28} & \multirow{2}{*}{$7.97e-1$}\\
   &  & {$1e-3$}& $2.08e-3$ &3 & $1.63e-1$ & &\\
 \midrule
\multirow{2}{*}{16384} & \multirow{2}{*}{$5.49e-7$} & {$1e-6$}& $2.73e-5$ &2 & $7.42e-1$ & \multirow{2}{*}{27} & \multirow{2}{*}{$4.68e+0$}\\
   &  & {$1e-3$}& $4.26e-3$ &3 & $9.02e-1$ & &\\
   \midrule
   \multirow{2}{*}{65536} & \multirow{2}{*}{$7.87e-7$} & {$1e-6$}& $9.56e-5$ &2 & $3.92e+0$ & \multirow{2}{*}{28} & \multirow{2}{*}{$2.72e+1$}\\
   &  & {$1e-3$}& $4.07e-3$ &3 & $4.83e+0$ & &\\
   \midrule
\multirow{2}{*}{262144} & \multirow{2}{*}{$8.90e-7$} & {$1e-6$}& $2.89e-4$ &2 & $2.07e+1$ & \multirow{2}{*}{27} & \multirow{2}{*}{$1.25e+2$}\\
   &  & {$1e-3$}& $5.57e-3$ &3 & $2.33e+1$ & &\\
\bottomrule
\end{tabular}
\caption{Numerical results for 1D uniform amplitude FIO using the approximate inverse $\Ga\Ka^*$ and the adjoint FIO matrix $\Ka^*$ as preconditioners for CG with tolerance $1e-8$.}
\label{1d_uniform}
\end{table}

   

\paragraph{Example 2. 1D FIO  with variable amplitude.} Here, we consider the 1D FIO from the previous example with variable periodic amplitude $a(x,\xi)$ in $x$ on the interval $[0,1]$, defined by a mixture of Gaussians of the form
\begin{align}
\sum_{k=0}^{n_k} e^{-\frac{(x-x_k)^2+(\xi-\xi_k)^2}{\sigma^2}},   
\end{align}
where $\sigma^{2}$ is set to $0.1$ or $0.05$.

\begin{figure}[h!]
\centering
 \includegraphics[width=0.75\textwidth,trim=0cm 0cm 0cm 1cm,clip]{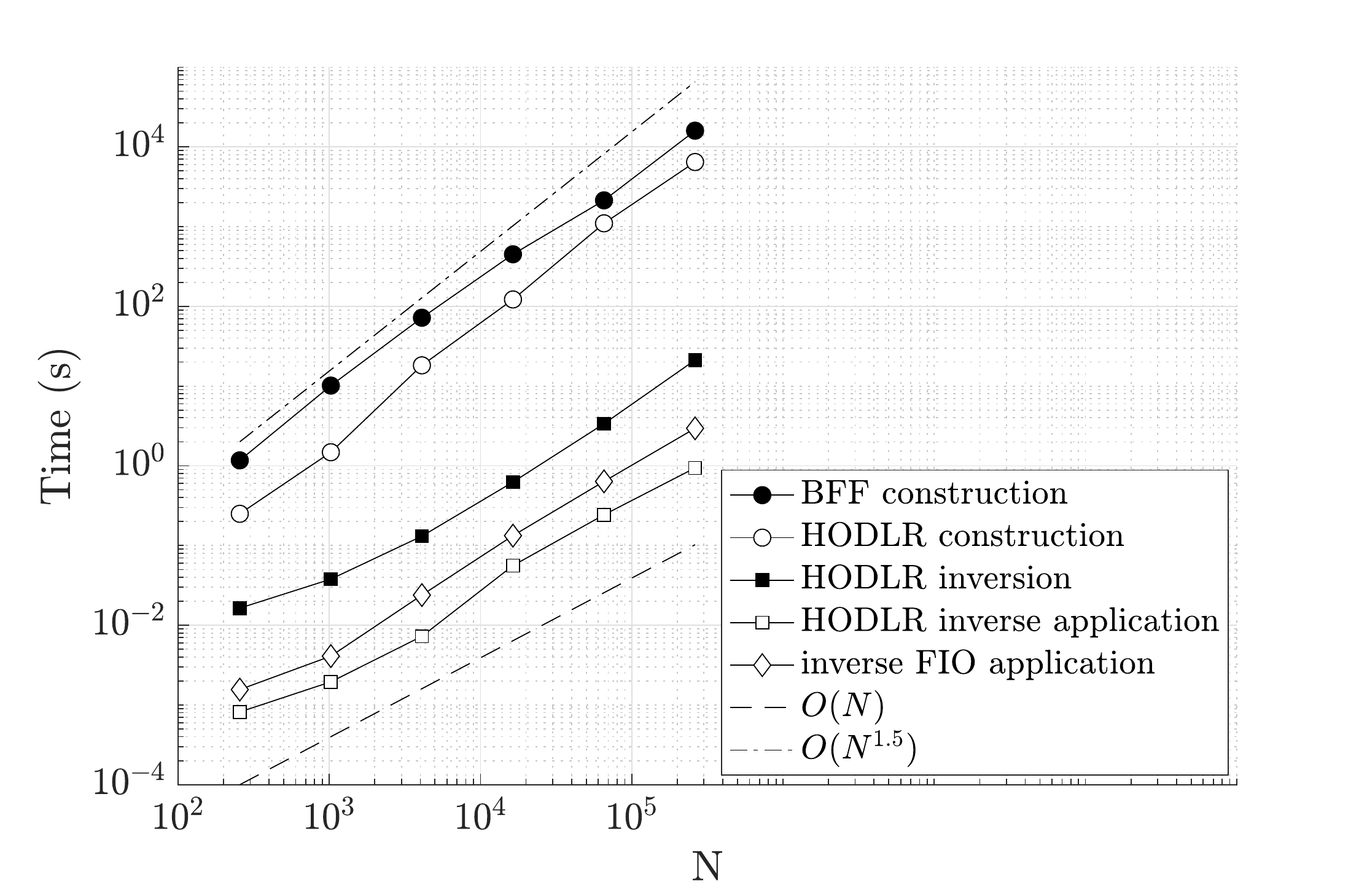}
\caption{Scaling results for 1D FIO from Example 2 with $\sigma^2=0.1$ and $\epsilon=0.001$. Timings are provided for each step of the algorithm to build the inverse factorization as well as the timings for matrix-vector multiplications with the resulting factorization. Reference $O(N)$ and $O(N^{3/2})$ lines are plotted, with dashed lines $--$ and $-\cdot-$, respectively.}
\label{1d_time}
\end{figure}

We report the results in Tables \ref{1d_variable} and \ref{1d_variable_2}. We observe that when
using a variable amplitude function, the adjoint FIO matrix $\Ka^*$ is not an efficient
preconditioner, since CG requires a very large number of iterations to converge, while the
approximate inverse $\Ga\Ka^*$ converges in just very few iterations, similar to Example
1. Additionally, the solve time with CG is 20 to 400 times lower depending on the tolerance
$\epsilon$ and the amplitude function $a(x,\xi)$ used. In particular, if $\sigma^2$ is set to $0.05$
as in Table \ref{1d_variable_2}, $\Ka^*$ becomes a very bad preconditioner, requiring the use of an
approximate factorization of $K^{-1}$.  We observe that constructing the approximate inverse
factorization takes $O(N^{3/2})$ complexity, while its application can be done in quasi-linear time,
as illustrated in Figure \ref{1d_time}.

\begin{table}[h!]
\centering
\begin{tabular}{c|c|cccc|cc}
\toprule
& $\Ka\approx \K $&\multicolumn{4}{|c|}{$\Ga\Ka^*\approx \K^{-1}$} & \multicolumn{2}{|c}{$\Ka^*\approx \Kt$}\\
\midrule
 $N$ & $e_a$ &$\epsilon$& $e_s$ & $n_i$ & $t_s$ & $n_i$ & $t_s$\\
\midrule
 \multirow{2}{*}{1024} & \multirow{2}{*}{$3.66e-7$} & {$1e-5$}& $6.13e-5$ &2 & $1.90e-2$ & \multirow{2}{*}{173} & \multirow{2}{*}{$6.81e-1$}\\
   &  & {$1e-3$}& $3.79e-2$ &3 & $2.71e-2$ & &\\
 \midrule
 \multirow{2}{*}{4096} & \multirow{2}{*}{$5.03e-7$} &{$1e-5$}& $1.38e-4$ &2 & $1.25e-1$ & \multirow{2}{*}{200} & \multirow{2}{*}{$5.80e+0$}\\
   &  & {$1e-3$}& $5.69e-3$ &3 & $1.73e-1$ & &\\
 \midrule
\multirow{2}{*}{16384} & \multirow{2}{*}{$6.69e-7$} & {$1e-5$}& $3.48e-4$ &2 & $7.38e-1$ & \multirow{2}{*}{205} & \multirow{2}{*}{$3.62e+1$}\\
   &  & {$1e-3$}& $8.52e-3$ &3 & $9.42e-1$ & &\\
   \midrule
   \multirow{2}{*}{65536} & \multirow{2}{*}{$8.49e-7$} & {$1e-5$}& $9.97e-4$ &3 & $5.18e+0$ & \multirow{2}{*}{207} & \multirow{2}{*}{$2.08e+2$}\\
   &  & {$1e-3$}& $8.45e-3$ &3 & $5.02e+0$ & &\\
   \midrule
\multirow{2}{*}{262144} & \multirow{2}{*}{$8.58e-7$} & {$1e-5$}& $3.92e-3$ & 3 &  $3.05e+1$ & \multirow{2}{*}{207} & \multirow{2}{*}{$1.10e+3$}\\
   &  & {$1e-3$}& $1.45e-2$ & 4& $3.61e+1$ & &\\
\bottomrule
\end{tabular}
\caption{Numerical results for 1D variable amplitude (with $\sigma^2 = 0.1$) FIO using the approximate inverse $\Ga\Ka^*$ and the adjoint FIO matrix $\Ka^*$ as preconditioners for CG with tolerance $1e-8$.}
\label{1d_variable}
\end{table}

\begin{table}[h!]
\centering
\begin{tabular}{c|c|cccc|cc}
  \toprule
& $\Ka\approx \K $&\multicolumn{4}{|c|}{$\Ga\Ka^*\approx \K^{-1}$} & \multicolumn{2}{|c}{$\Ka^*\approx \Kt$}\\
\midrule
 $N$ & $e_a$ &$\epsilon$& $e_s$ & $n_i$ & $t_s$ & $n_i$ & $t_s$\\
\midrule
 \multirow{2}{*}{1024} & \multirow{2}{*}{$3.33e-7$} & {$1e-4$}& $3.11e-3$ &3 & $3.06e-2$ & \multirow{2}{*}{1020} & \multirow{2}{*}{$4.02e+0$}\\
   &  & {$1e-3$}& $1.02e-2$ &4 & $3.84e-2$ & &\\
 \midrule
 \multirow{2}{*}{4096} & \multirow{2}{*}{$4.10e-7$} &{$1e-4$}& $5.57e-3$ &3 & $1.88e-1$ & \multirow{2}{*}{2221} & \multirow{2}{*}{$6.43e+1$}\\
   &  & {$1e-3$}& $1.24e-1$ &4 & $2.22e-1$ & &\\
 \midrule
\multirow{2}{*}{16384} & \multirow{2}{*}{$7.08e-7$} & {$1e-4$}& $1.37e-2$ &4 & $1.20e+0$ & \multirow{2}{*}{2517} & \multirow{2}{*}{$4.29e+2$}\\
   &  & {$1e-3$}& $2.45e-2$ &4 & $1.33e+0$ & &\\
   \midrule
   \multirow{2}{*}{65536} & \multirow{2}{*}{$7.51e-7$} & {$1e-4$}& $1.87e-2$ &4 & $6.15e+0$ & \multirow{2}{*}{2630} & \multirow{2}{*}{$2.34e+3$}\\
   &  & {$1e-3$}& $4.46e-2$ &5 & $7.98e+0$ & &\\
   \midrule
   \multirow{2}{*}{262144} & \multirow{2}{*}{$7.29e-7$} & {$1e-4$}& $6.36e-2$ &5 & $3.38e+1$ & \multirow{2}{*}{2668} & \multirow{2}{*}{$1.13e+4$}\\
   &  & {$1e-3$}& $1.75e-1$ &7 & $4.39e+1$ & &\\
   \midrule
\end{tabular}
\caption{Numerical results for 1D variable amplitude (with $\sigma^2 = 0.05$) FIO using the approximate inverse $\Ga\Ka^*$ and the adjoint FIO matrix $\Ka^*$ as preconditioners for CG with tolerance $1e-8$.}
\label{1d_variable_2}
\end{table}

\paragraph{Example 3. 2D generalized Radon transform.} Finally, we consider the 2D analog of the numerical Example 1, 
\begin{align}
  \label{eq_2d_FIO}
  u(x) = \sum_{\xi\in \Omega} e^{2\pi i \Phi(x,\xi)}\hat{f}(\xi), \quad x\in X
\end{align} 
with phase function $\Phi(x,\xi)$ defined as
\begin{align}
\label{eq_2d_FIO_fun}
\begin{split}
\Phi(x,\xi) = x\cdot\xi + \sqrt{c_1^{2}(x)\xi_1^2 + c_2^{2}(x)\xi_2^2},\\
c_1(x) = (2 + \sin(2\pi x_1)\sin(2\pi x_2))/16, \\
c_2(x) = (2 + \cos(2\pi x_1)\cos(2\pi x_2))/16,
\end{split}
\end{align}
and
\begin{align}
  \label{eq_2d_Omega}
  X = \left\{x=\Big(\frac{n_1}{n}, \frac{n_2}{n}\Big),0\leq n_1,n_2< n,\text{ with }n_1,n_2\in\mathbb{Z}\right\},\\
  \Omega = \left\{\xi=(n_1,n_2), -\frac{n}{2}\leq n_1,n_2< \frac{n}{2},\text{ with }n_1,n_2\in\mathbb{Z}\right\},
\end{align} 
with $n$ being the number of points in each dimension and $N = n^2$. This leads again to the system $u = Kf$ with $K = (e^{2\pi i \Phi(x,\xi)})_{x\in X,\xi\in\Omega}$, for which we are interested in the inversion of matrix $K$ and the solution $f = K^{-1}u$. This FIO models an integration over ellipses centered at a point $x\in X$ with axis lengths defined by $c_1(x)$ and $c_2(x)$.

\begin{figure}[h!]
\centering
 \includegraphics[width=0.75\textwidth,trim=0cm 0cm 0cm 0.7cm,clip]{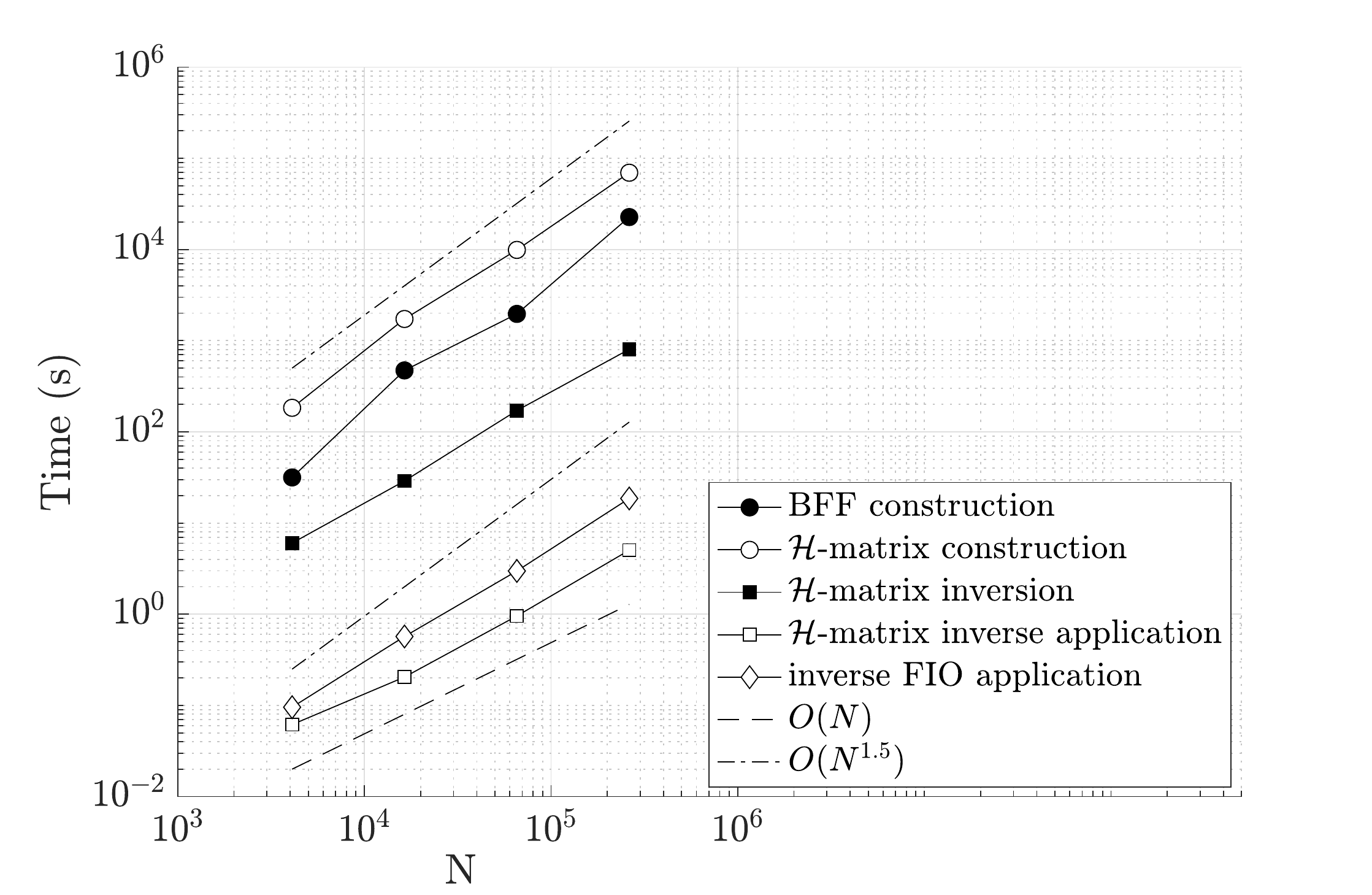}
\caption{Scaling results for 2D FIO with $\epsilon=0.001$. Timings are provided for constructing BFF
  and $\mathcal{H}$-matrix as well as for inverting the $\mathcal{H}$-matrix and performing
  matrix-vector multiplications with the factorizations of $\Ga$ and $\Ga\Ka^*$. Reference $O(N)$ and
  $O(N^{3/2})$ lines are plotted, with dashed lines $--$ and $-\cdot-$ respectively.}
\label{2d_time}
\end{figure}

One can use the butterfly factorization to obtain a factorization $\Ka$ of $K$ with forward error
$e_a$ between $10^{-4}$ and $10^{-5}$. As observed in Table \ref{2d_FIO}, using tolerance
$\epsilon=10^{-3}$ for the peeling algorithm and the $\mathcal{H}$-matrix inversion we can find an
approximate inverse factorization that is $\epsilon$-accurate. This factorization of the inverse can
be used as direct solver or as a preconditioner for CG, which converges in just 3 iterations, as
opposed to the 20 iterations needed when only $\Ka^*$ is used.  Notice also that the total timing of
CG is 4 times lower.  The timings for the different parts of the algorithm are illustrated in Figure
\ref{2d_time}. We observe $O(N^{3/2})$ time complexity for building the factorization and almost
$O(N)$ time for performing matrix–vector multiplication. The algorithm proposed for
$\mathcal{H}$-matrix inversion experimentally shows quasi-linear complexity.

\begin{table}[h!]
\centering
\begin{tabular}{c|c|cccc|cc}
\toprule
& $\Ka\approx \K $&\multicolumn{4}{|c|}{$\Ga\Ka^*\approx \K^{-1}$} & \multicolumn{2}{|c}{$\Ka^*\approx \Kt$}\\
\midrule
 $N$ & $e_a$ &$\epsilon$& $e_s$ & $n_i$ & $t_s$ & $n_i$ & $t_s$\\
\midrule
$64^2$ & $9.62e-5$ & $1e-3$ & $2.13e-3$ &3 & $7.09e-1$ & 20 & $1.36e+0$\\
$128^2$ & $9.72e-5$ & $1e-3$ & $2.64e-3$ &3 & $4.25e+0$ & 21 & $1.56e+1$\\
$256^2$ & $2.05e-4$ & $1e-3$ & $4.75e-3$ &3 & $2.27e+1$ & 22 & $7.6e+1$\\
$512^2$ & $2.27e-4$ & $1e-3$ & $9.46e-3$ &4 & $ 1.87e+2$ & 22 & $6.5e+2$\\
\bottomrule
\end{tabular}
\caption{Numerical results for 2D FIO using the approximate inverse $\Ga\Ka^*$ and the adjoint FIO matrix $\Ka^*$ as preconditioners with tolerance $1e-8$.}
\label{2d_FIO}
\end{table}

\section{Conclusions}\label{sec:conclusions}

In this paper, we have introduced an algorithm to obtain a factorization of inverse discrete FIO
matrices $K\in\mathbb{C}^{N\times N}$ in 1D and 2D that can be computed in $O(N^{3/2})$ time and
space complexity. The algorithm combines the butterfly factorization with hierarchical matrices, to
exploit the complimentary low-rank condition of discrete FIO matrices $K$ and the numerical low-rank
subblocks that appear in the corresponding hermitian matrix $K^*K$.  With such factorization, the
inverse matrix $K^{-1}$ is represented as a product of $O(\log N )$ sparse matrices, and therefore
it can be applied in quasi-linear time, $O(N\log N)$.  This factorization can be used either as a
direct solver or as a preconditioner for iterative methods depending on the accuracy required.  When
used as a preconditioner the factorization can reduce the number of iterations of iterative methods
such as CG or GMRES, which otherwise can take $O(N)$ iterations to converge. We have also proposed a new algorithm to invert HODLR and $\mathcal{H}$-matrix representations based on recursive skeletonization and HIF, that construct approximate inverse factorizations in quasi-linear time.

Notice that the three main steps of the algorithm (butterfly factorization, $\mathcal{H}$-matrix construction and $\mathcal{H}$-matrix inversion) can be easily parallelized. The butterfly factorization and $\mathcal{H}$-matrix inversion algorithms are organized following a tree structure, where at a each level every node can be processed independently of the rest, as done in \cite{DHIF}. Similarly, the $\mathcal{H}$-matrix construction for the 2D case can be easily parallelized at each level of the tree with 64 processors, one for each subset $\mathcal{P}_{p,q}$ in \eqref{eq_spq}.

The algorithm proposed here can also be used on higher dimensions ($d>2$) by partitioning the domain in $d-$dimensional hypercubes, as far as admissible blocks on $K^*K$ can be considered low-rank.

\section*{Acknowledgments}\label{acknowlegements}
The work of J.F. is supported by the Stanford Graduate Fellowship. The work of L.Y. is partially
supported by U.S. Department of Energy, Office of Science, Office of Advanced Scientific Computing
Research, Scientific Discovery through Advanced Computing (SciDAC) program and the National Science
Foundation under award DMS-1818449. We thank Abeynaya Gnanasekaran for valuable conversations on the
numerical algorithms.

\end{document}